\documentclass[12pt]{article}
\usepackage{amsfonts}
\usepackage{full page,
amssymb, amscd, amsmath,graphicx, 
}

 \title{{\bf Associative algebras for (logarithmic) twisted modules for a vertex operator algebra}}
 \author{Yi-Zhi Huang and Jinwei Yang}
    \date{}
    \begin{document}
    \bibliographystyle{alpha}
    \maketitle
\newtheorem{thm}{Theorem}[section]
\newtheorem{defn}[thm]{Definition}
\newtheorem{prop}[thm]{Proposition}
\newtheorem{cor}[thm]{Corollary}
\newtheorem{lemma}[thm]{Lemma}
\newtheorem{rema}[thm]{Remark}
\newtheorem{app}[thm]{Application}
\newtheorem{prob}[thm]{Problem}
\newtheorem{conv}[thm]{Convention}
\newtheorem{conj}[thm]{Conjecture}
\newtheorem{cond}[thm]{Condition}
    \newtheorem{exam}[thm]{Example}
\newtheorem{assum}[thm]{Assumption}
     \newtheorem{nota}[thm]{Notation}
\newcommand{\halmos}{\rule{1ex}{1.4ex}}
\newcommand{\pfbox}{\hspace*{\fill}\mbox{$\halmos$}}

\newcommand{\nn}{\nonumber \\}

 \newcommand{\res}{\mbox{\rm Res}}
 \newcommand{\ord}{\mbox{\rm ord}}
\renewcommand{\hom}{\mbox{\rm Hom}}
\newcommand{\edo}{\mbox{\rm End}\ }
 \newcommand{\pf}{{\it Proof.}\hspace{2ex}}
 \newcommand{\epf}{\hspace*{\fill}\mbox{$\halmos$}}
 \newcommand{\epfv}{\hspace*{\fill}\mbox{$\halmos$}\vspace{1em}}
 \newcommand{\epfe}{\hspace{2em}\halmos}
\newcommand{\nord}{\mbox{\scriptsize ${\circ\atop\circ}$}}
\newcommand{\wt}{\mbox{\rm wt}\ }
\newcommand{\swt}{\mbox{\rm {\scriptsize wt}}\ }
\newcommand{\lwt}{\mbox{\rm wt}^{L}\;}
\newcommand{\rwt}{\mbox{\rm wt}^{R}\;}
\newcommand{\slwt}{\mbox{\rm {\scriptsize wt}}^{L}\,}
\newcommand{\srwt}{\mbox{\rm {\scriptsize wt}}^{R}\,}
\newcommand{\clr}{\mbox{\rm clr}\ }
\newcommand{\tr}{\mbox{\rm Tr}}
\newcommand{\A}{\mathbf{A}}
\renewcommand{\L}{\mathbf{L}}
\newcommand{\C}{\mathbb{C}}
\newcommand{\Z}{\mathbb{Z}}
\newcommand{\R}{\mathbb{R}}
\newcommand{\Q}{\mathbb{Q}}
\newcommand{\N}{\mathbb{N}}
\newcommand{\CN}{\mathcal{N}}
\newcommand{\F}{\mathcal{F}}
\newcommand{\I}{\mathcal{I}}
\newcommand{\V}{\mathcal{V}}
\newcommand{\one}{\mathbf{1}}
\newcommand{\BY}{\mathbb{Y}}
\newcommand{\ds}{\displaystyle}

        \newcommand{\ba}{\begin{array}}
        \newcommand{\ea}{\end{array}}
        \newcommand{\be}{\begin{equation}}
        \newcommand{\ee}{\end{equation}}
        \newcommand{\bea}{\begin{eqnarray}}
        \newcommand{\eea}{\end{eqnarray}}
         \newcommand{\lbar}{\bigg\vert}
        \newcommand{\p}{\partial}
        \newcommand{\dps}{\displaystyle}
        \newcommand{\bra}{\langle}
        \newcommand{\ket}{\rangle}

        \newcommand{\ob}{{\rm ob}\,}
        \renewcommand{\hom}{{\rm Hom}}

\newcommand{\Y}{\mathcal{Y}}

\newcommand{\dlt}[3]{#1 ^{-1}\delta \bigg( \frac{#2 #3 }{#1 }\bigg) }

\newcommand{\dlti}[3]{#1 \delta \bigg( \frac{#2 #3 }{#1 ^{-1}}\bigg) }

 \makeatletter
\newlength{\@pxlwd} \newlength{\@rulewd} \newlength{\@pxlht}
\catcode`.=\active \catcode`B=\active \catcode`:=\active
\catcode`|=\active
\def\sprite#1(#2,#3)[#4,#5]{
   \edef\@sprbox{\expandafter\@cdr\string#1\@nil @box}
   \expandafter\newsavebox\csname\@sprbox\endcsname
   \edef#1{\expandafter\usebox\csname\@sprbox\endcsname}
   \expandafter\setbox\csname\@sprbox\endcsname =\hbox\bgroup
   \vbox\bgroup
  \catcode`.=\active\catcode`B=\active\catcode`:=\active\catcode`|=\active
      \@pxlwd=#4 \divide\@pxlwd by #3 \@rulewd=\@pxlwd
      \@pxlht=#5 \divide\@pxlht by #2
      \def .{\hskip \@pxlwd \ignorespaces}
      \def B{\@ifnextchar B{\advance\@rulewd by \@pxlwd}{\vrule
         height \@pxlht width \@rulewd depth 0 pt \@rulewd=\@pxlwd}}
      \def :{\hbox\bgroup\vrule height \@pxlht width 0pt depth
0pt\ignorespaces}
      \def |{\vrule height \@pxlht width 0pt depth 0pt\egroup
         \prevdepth= -1000 pt}
   }
\def\endsprite{\egroup\egroup}
\catcode`.=12 \catcode`B=11 \catcode`:=12 \catcode`|=12\relax
\makeatother

\def\hboxtr{\FormOfHboxtr} 
\sprite{\FormOfHboxtr}(25,25)[0.5 em, 1.2 ex] 

:BBBBBBBBBBBBBBBBBBBBBBBBB | :BB......................B |
:B.B.....................B | :B..B....................B |
:B...B...................B | :B....B..................B |
:B.....B.................B | :B......B................B |
:B.......B...............B | :B........B..............B |
:B.........B.............B | :B..........B............B |
:B...........B...........B | :B............B..........B |
:B.............B.........B | :B..............B........B |
:B...............B.......B | :B................B......B |
:B.................B.....B | :B..................B....B |
:B...................B...B | :B....................B..B |
:B.....................B.B | :B......................BB |
:BBBBBBBBBBBBBBBBBBBBBBBBB |

\endsprite
\def\shboxtr{\FormOfShboxtr} 
\sprite{\FormOfShboxtr}(25,25)[0.3 em, 0.72 ex] 

:BBBBBBBBBBBBBBBBBBBBBBBBB | :BB......................B |
:B.B.....................B | :B..B....................B |
:B...B...................B | :B....B..................B |
:B.....B.................B | :B......B................B |
:B.......B...............B | :B........B..............B |
:B.........B.............B | :B..........B............B |
:B...........B...........B | :B............B..........B |
:B.............B.........B | :B..............B........B |
:B...............B.......B | :B................B......B |
:B.................B.....B | :B..................B....B |
:B...................B...B | :B....................B..B |
:B.....................B.B | :B......................BB |
:BBBBBBBBBBBBBBBBBBBBBBBBB |

\endsprite

\date{}
\maketitle

\begin{abstract}
We construct two associative algebras from a vertex operator algebra $V$ and a 
general automorphism $g$ of $V$. The first, called $g$-twisted zero-mode algebra, is a 
subquotient of what we call $g$-twisted universal enveloping algebra of $V$. These
algebras are generalizations of the corresponding algebras introduced and studied by 
Frenkel-Zhu and Nagatomo-Tsuchiya
in the (untwisted) case that $g$ is the identity. The other 
is a generalization of the $g$-twisted version of Zhu's algebra for suitable $g$-twisted 
modules constructed by Dong-Li-Mason when the order of $g$ is finite. 
We are mainly interested in $g$-twisted $V$-modules introduced by the first author 
in the case that $g$ is of infinite order and does not act 
on $V$ semisimply. In this case, twisted vertex operators in general 
involve the logarithm of the variable. 
We construct  functors between categories of suitable modules 
for these associative algebras and categories of suitable (logarithmic) $g$-twisted 
$V$-modules. Using these functors, we prove that the
$g$-twisted zero-mode algebra and the $g$-twisted generalization of Zhu's algebra
are in fact isomorphic. 
\end{abstract}



\renewcommand{\theequation}{\thesection.\arabic{equation}}
\renewcommand{\thethm}{\thesection.\arabic{thm}}
\setcounter{equation}{0} \setcounter{thm}{0} 

\section{Introduction}

Orbifold conformal field theories play an important role in mathematics and physics. They 
are examples of conformal field theories constructed from known conformal field theories and 
automorphisms of these known ones. 
Mathematically, the study of orbifold conformal field theories 
can be reduced to the study of the twisted
representation theory of vertex operator algebras, that is, 
the theory of representations of vertex operator algebras twisted by automorphisms.
In fact, the first example of orbifold conformal field theories is the conformal field 
theory corresponding to the moonshine module vertex operator algebra constructed by 
Frenkel, Lepowsky and Meurman \cite{FLM}. 
One conjecture proposed by the first author in this twisted representation theory of vertex operator algebras
is that $g$-twisted modules for a vertex operator algebra $V$ satisfying suitable conditions and 
elements $g$ of a group $G$ of automorphisms of $V$ equipped with 
twisted (logarithmic) intertwining operators among $g$-twisted modules for $g\in G$ have a structure of $G$-equivariant 
intertwining operator algebras satisfying a modular invariance property. This conjecture in fact 
implies that the category of $g$-twisted $V$-modules for $g\in G$ has a structure of $G$-crossed (tensor) category in the sense of
Turaev \cite{T}.

Automorphisms of a vertex operator algebra can be of finite or infinite orders. 
An element of the Monster group gives an automorphism of the moonshine module 
vertex operator algebra \cite{FLM}. Such an automorphism 
is of finite order. An element of a simply connected finite-dimensional Lie group gives an automorphism 
of the vertex operator algebra associated to the affine Lie algebra of the Lie algebra of the Lie group \cite{FLM} \cite{FZ}.
Such an automorphism is in general of infinite order. Moreover, in general  it might not act on the vertex 
operator algebra semisimply. 

Twisted modules for a vertex operator algebra associated to an automorphism of finite orders
were introduced in the construction of 
the moonshine module vertex operator algebra in
\cite{FLM}. In \cite{H}, 
the first author introduced a notion of twisted module associated to 
a general automorphism  whose order does not have to be finite.
One important new feature of twisted modules
in the general case is that twisted vertex operators might involve the logarithm of the variable (logarithmic 
$g$-twisted module). 
For a historical discussion of twisted modules
for vertex operator algebras and the relevant references, see \cite{H}. See also \cite{B} for the Jacobi 
identity for a suitable component of the twisted vertex operators that we shall derive from the axioms in 
\cite{H} and use in this paper. 

In the representation theory of Lie algebras, one of the most important 
tools is the universal enveloping algebra of a Lie algebra. The representation theory of a
Lie algebra is equivalent to the representation theory of its universal enveloping algebra. 
It is natural to try to find associative algebras that play the similar roles in the representation theory 
of vertex operator algebras. 

Let $V$ be a vertex operator algebra. 
Frenkel and Zhu in \cite{FZ} first constructed the universal enveloping algebra $U(V)$ of $V$. 
$V$-modules of the weakest type  (even more general than the so called weak $V$ modules) are equivalent to $U(V)$-modules. 
Though $U(V)$  is very natural, it is not very useful because in the representation theory of vertex  operator algebras,
we are interested mostly in $V$-modules with suitable lower bounded 
gradings. In \cite{Z}, Zhu constructed an associative algebra $A(V)$ on a quotient of $V$ with 
the product obtained from the vertex operator map. He also constructed implicitly functors 
between the category of $V$-modules with suitable lower bounded gradings and the category of $A(V)$-modules. 
In \cite{H1}, motivated by the geometry of vertex operator algebras, 
the first author constructed an associative algebra $\tilde{A}(V)$ using a product and a 
quotient that look very different from those defining $A(V)$. It was also proved in \cite{H1} that $\tilde{A}(V)$
is in fact isomorphic to $A(V)$ by constructing  explicitly an isomorphism between them using 
the conformal transformation $\frac{1}{2\pi i}\log (1+2\pi i z)$. 

To study twisted modules for vertex operator algebras, we would also like to 
find suitable associative algebras such that the study of twisted modules can be reduced to the study of 
modules for these associative algebras. In the case that the automorphism $g$ of the vertex 
operator algebra $V$ is of finite order,   Dong, Li and Mason in \cite{DLM} generalized Zhu's construction 
to obtain an associative algebra $A_{g}(V)$ 
together with functors between the categories of suitable $g$-twisted $V$-modules and $A_{g}(V)$-modules. 
This associative algebra $A_{g}(V)$  is in fact a subalgebra of Zhu's algebra for the fixed point vertex operator 
algebra under $g$.  The construction 
in  \cite{DLM}  cannot be generalized directly to the case that the order of  $g$  is infinite. 

In the present paper,  for a vertex operator algebra $V$ and an automorphism $g$ of $V$, we first construct 
a $g$-twisted universal enveloping algebra $U_{g}(V)$ such that $g$-twisted $V$-modules of the weakest type 
 (even more general than $g$-twisted weak $V$-modules) are equivalent to $U_{g}(V)$-modules. 
The construction is a straightforward generalization of the construction of Frenkel-Zhu to the twisted 
case.  We take the tensor algebra of the affinization of $V$ with suitable powers of the variable and take a topological completion
so that suitable infinite sums are allowed, and then divide this algebra by all the relations corresponding to the identities 
that should hold for any type of $g$-twisted modules. 

Just as in the case of $V$-modules, $U_{g}(V)$  is not very 
useful because we are interested mostly in $g$-twisted $V$-modules with suitable lower bounded 
gradings. For $g$-twisted $V$-modules graded by conformal weights, there exist lowest weight spaces.
In the case that such a $g$-twisted $V$-module is irreducible, the lowest weight space determines the module completely. 
It is natural to expect that there exists an associative algebra obtained from $U_{g}(V)$ such that 
the lowest weight space of an irreducible  $g$-twisted $V$-module is a module for this algebra. 

In this paper, 
using $U_{g}(V)$, we construct an associative algebra $Z_{g}(V)$ satisfying 
this property and call it $g$-twisted zero-mode algebra,  although  the more appropriate name
for this algebra is probably $g$-twisted imaginary-mode algebra or $g$-twisted zero-real-mode algebra. 
In fact, $U_{g}(V)$ is graded by conformal weights. Take the 
subalgebra of $U_{g}(V)$ generated by elements of imaginary weights 
and then take the quotient of this algebra by the ideal generated by 
elements of the form $uv$ where $u, v\in U_{g}(V)$ have weights $k+m$ and $-m$ for 
$k\in \mathbb{I}$ (the set of imaginary numbers) and $m\in \C$ satisfying $\Re(m)>0$. 
This quotient algebra is our $g$-twisted zero-mode algebra $Z_{g}(V)$.

We also construct  an associative algebra $\tilde{A}_{g}(V)$, generalizing the algebra $\tilde{A}(V)$ in \cite{H1},
and an isomorphic algebra $A_{g}(V)$, generalizing the algebra constructed in \cite{DLM} for $g$ of finite order. 
In fact, in this case, $\tilde{A}_{g}(V)$ is more natural to work with. 
We give our constructions and proofs mostly for $\tilde{A}_{g}(V)$. The construction and results for $A_{g}(V)$
can be easily derived from those for $\tilde{A}_{g}(V)$ using the isomorphism between them. 
One interesting feature in the case that the order of $g$ is infinite is that $A_{g}(V)$
 is in general not a subalgebra 
of Zhu's algebra for the fixed point vertex operator subalgebra under $g$. It is instead a subalgebra 
of Zhu's algebra for the fixed point vertex operator subalgebra under the semisimple part of $g$. 

Our main results in this paper are the constructions of functors between the categories of suitable 
$g$-twisted $V$-modules, suitable $Z_{g}(V)$-modules and suitable $\tilde{A}_{g}(V)$-modules (or $A_{g}(V)$-modules). 
As in the case that the order of $g$ is finite, a suitable $g$-twisted $V$-module $W$ has a subspace $\Omega_{g}(W)$
on which the components of the twisted vertex operators whose weights have negative real parts  act as $0$. 
It is easy to prove from the definitions that $\Omega_{g}(W)$ is a $Z_{g}(V)$-,  $\tilde{A}_{g}(V)$- or $A_{g}(V)$-module.
We then construct right inverses of these functors and derive some consequences that will be useful for future study.
Finally, using these results, we prove that $Z_{g}(V)$, $\tilde{A}_{g}(V)$ and $A_{g}(V)$ are isomorphic to each other.

In the case that $g$ is the identity, 
the $g$-twisted zero-mode algebra becomes the zero-mode algebra, which appeared first in physics (see \cite{BN})
and introduced rigorously in \cite{FZ}. In \cite{NT},  Nagatomo and Tsuchiya \cite{NT}
proved that the zero-mode algebra  is isomorphic to Zhu's algebra. As a special case (the case that $g$ is the identity),
our proof that $Z_{g}(V)$ and $A_{g}(V)$ are isomorphic gives in particular a different proof that the 
zero-mode algebra  is isomorphic to Zhu's algebra.

It seems to be easier to calculate $Z_{g}(V)$ than to calculate $A_{g}(V)$ or $\tilde{A}_{g}(V)$. On the other hand,
$\tilde{A}_{g}(V)$ is more natural for the study of modular invariance. $A_{g}(V)$ is in some sense a bridge 
connecting them.  We expect that $Z_{g}(V)$ will be especially useful for the construction and
study of twisted modules and twisted intertwining operators. 

There are also higher zero-mode algebras $Z_{g, n}(V)$ for $n\in \N$ as the quotient of the
subalgebra of $U_{g}(V)$ generated by elements of imaginary weights by the ideal generated by 
elements of the form $uv$ where $u, v\in U_{g}(V)$ have weights $k+m$ and $-m$ for $k\in\mathbb{I}$ and
$m\in \C$ satisfying $\Re(m)>n$. 
The $g$-twisted zero-mode algebra $Z_{g}$ discussed above is in fact $Z_{g, 0}(V)$. 
We can also generalize higher Zhu's algebras $A_{n}(V)$ introduced by Dong, Li and Mason in \cite{DLM2} to 
associative algebras $A_{g, n}(V)$ for $n\in \N$ such that $A_{g, 0}(V)=A_{g}(V)$ and $A_{1, n}(V)=A_{n}(V)$
when $g$ is the identity $1$. There are also $\tilde{A}_{g, n}(V)$  for $n\in \N$ generalizing $\tilde{A}_{g}(V)$. 
These algebras will be important for the study of twisted logarithmic intertwining operators. 
But for simplicity, we shall study only $Z_{g}(V)$, $\tilde{A}_{g}(V)$ and $A_{g}(V)$ in this paper. 
The general case will be discussed in another paper.

The definition of $g$-twisted module in \cite{H} was formulated using the duality properties of the twisted vertex operators.
It should be possible to use the duality properties and immediate consequences to construct the associative algebras
and functors above (see \cite{HY} for a construction of the functors between the category of modules for Zhu's algebra
$A(V)$ and the category of  suitable $V$-modules without using the commutator formula for
modules).
On the other hand, since twisted vertex operators are special (logarithmic) intertwining operators, suitable components 
of twisted vertex operators might satisfy some Jacobi-type identity as is discussed for intertwining operators by the first author 
in \cite{H0}. A Jacobi-type identity will give  formulas that can be used to simplify our constructions. 
Indeed, a Jacobi identity 
 for suitable components of twisted vertex operators was obtained 
by Bakalov in \cite{B}. Bakalov gave in \cite{B} a different definition of $g$-twisted $V$-module based on an 
associator formula for the twisted vertex operators. From this definition, he derived a Jacobi identity 
for the map obtained by taking the non-logarithmic component of the twisted vertex operator map, that is, 
the constant term in the twisted vertex operator map when we view the twisted vertex operator map as a 
polynomial in the logarithm of the variable.
It was stated in \cite{B} that a $g$-twisted $V$-module satisfying the original definition given by the first author in \cite{H} 
indeed satisfies the definition given in \cite{B}. But no proof was given there. 

Note that the twisted vertex operators given in the definition in \cite{H} are the objects that
we are interested in orbifold conformal field theory, not their non-logarithmic components. 
For example, in the future, we would like to prove the 
modular invariance for the space of $q$-(pseudo-)traces of twisted vertex operators for $g$-twisted $V$-modules, 
not for the space of $q$-(pseudo-)traces of non-logarithmic components of twisted vertex operators for 
$g$-twisted $V$-modules. Thus before we can use the Jacobi identity in \cite{B}, we need to 
first show that for suitable $g$-twisted $V$-modules, these two definitions are equivalent. We give a proof of this equivalence in 
this paper. We then use freely the formulas derived from this Jacobi identity to give our constructions and proofs.

The present paper is organized as follows: In Section 2, we recall several variants of the notion of 
$g$-twisted $V$-module. We derive Bakalov's  Jacobi identity for the non-logarithmic components of twisted 
vertex operators and some other useful properties in this section. In this section, we also prove that 
when the other conditions for (generalized) $g$-twisted $V$-modules hold, the duality property for the 
twisted vertex operators given in \cite{H} is equivalent to the Jacobi identity for non-logarithmic components
of twisted vertex operators given in \cite{B}. We give the constructions of $g$-twisted universal enveloping algebra
$U_{g}(V)$ and the zero-mode algebra $Z_{g}(V)$ in Section 3. The constructions of $\tilde{A}_{g}(V)$ and $A_{g}(V)$ 
are given in Section 4. Our main results, the constructions of functors between different categories, are given in 
Section 5. In Section 6, we give the proof that $Z_{g}(V)$, $A_{g}(V)$ and $\tilde{A}_{g}(V)$ are isomorphic to 
each other.

\paragraph{Acknowledgments}
We would like to thank Robert McRae for comments on the formulation of the duality property in the 
definition of twisted modules.

\renewcommand{\theequation}{\thesection.\arabic{equation}}
\renewcommand{\thethm}{\thesection.\arabic{thm}}
\setcounter{equation}{0} \setcounter{thm}{0} 

\section{Twisted modules}

In this paper, we fix a vertex operator algebra $(V, Y, {\bf 1}, \omega)$ and an automorphism  $g$
of $V$. We shall use the convention that $\log z=\log |z|+\arg z \sqrt{-1}$
where $0\le \arg z<2\pi$. We shall also use $l_{p}(z)$ to denote $\log z+2 \pi p\sqrt{-1}$ for $p\in \Z$. 
Besides the usual notations $\C$, $\R$, $\Z$, $\Z_{+}$ and $\N$ for the sets of complex numbers, 
real numbers, integers, positive integers and natural numbers, respectively, we also use $\mathbb{I}$ and 
$\overline{\C}_{+}$ to denote the set of imaginary numbers and the closed right half complex plane. 

In this section, we recall the definition of several variants of the notion of 
$g$-twisted $V$-module and prove a number of useful properties for these $g$-twisted $V$-modules. 
Many of the formulas proved in this section were obtained first by Bakalov \cite{B} from a different 
definition of $g$-twisted $V$-module. Here we prove them starting from the original definition in \cite{H}. 
We also prove that 
when the other conditions for suitable $g$-twisted $V$-modules hold, the duality property for the 
twisted vertex operators given in \cite{H} is equivalent to the Jacobi identity for non-logarithmic components
of twisted vertex operators given in \cite{B}. 

For $\alpha \in \C/\Z$, we denote by
$$
V^{[\alpha]} = \{u \in V\;|\; (g-e^{2\pi\sqrt{-1}\alpha})^{\Lambda} u =0 \;\text{for some}\; \Lambda\in \Z_{+}\}
$$
the generalized eigenspace with eigenvalue $e^{2\pi \sqrt{-1}\alpha}$ for $g$. 
For $\alpha\in \C/\Z$, 
there is a unique $a\in [0, 1)+\mathbb{I}$ such that $a+\Z=\alpha$. 
We call $a\in [0, 1)+\mathbb{I}$ a
{\it $g$-weight of $V$} if $V^{[a+\Z]} \neq 0$. 
We use $P(V)$ to denote the
set of all the $g$-weights of $V$.  Then
$V=\coprod_{n\in \C, a\in P(V)}V_{(n)}^{[a+\Z]}$. Let $\A_{V}: V\to V$ be the linear map defined 
by $\A_{V} v=av$ for $v\in V^{[a+\Z]}$ and $a\in P(V)$. 

\begin{defn}\label{d}
{\rm  A {\it  $\overline{\C}_{+}$-graded weak $g$-twisted
$V$-module} is a $\overline{\C}_{+}\times \mathbb{C}/\mathbb{Z}$-graded
vector space $W = \coprod_{n \in \overline{\C}_{+}, \alpha \in
\mathbb{C}/\mathbb{Z}} W_{n}^{[\alpha]}$ (graded by {\it $\overline{\C}_{+}$-degrees} and {\it $g$-weights}) 
equipped with a linear
map
\begin{align*}
Y^g: &V\otimes W \rightarrow W\{x\}[{\rm log} x],\\
&v \otimes w \mapsto Y^g(v, x)w
\end{align*}
 and an action of $g$ satisfying the following conditions:
\begin{enumerate}

\item The {\it equivariance property}: For $p \in \mathbb{Z}$, $z
\in \mathbb{C}^{\times}$, $v \in V$ and $w \in W$, $Y^{g; p + 1}(gv,
z)w = Y^{g; p}(v, z)w$, where for $p \in \mathbb{Z}$, 
$$Y^{g; p}(v, z)w=Y^{g}(v, x)w\lbar_{x^{n}=e^{nl_{p}(z)},\; \log x=l_{p}(z)}$$
is the $p$-th analytic branch of $Y^g$.

\item The {\it identity property}: For $w \in W$, $Y^g({\bf 1}, x)w
= w$.

\item  The {\it duality property}: Let $W^{'} = \coprod_{n \in \overline{\C}_{+}, \alpha \in
\mathbb{C}/\mathbb{Z}} \left(W_{n}^{[\alpha]}\right)^{\ast}$ and, for $n \in
\mathbb{C}$, $\pi_n : W \rightarrow W_{n}=\coprod_{\alpha\in \C/\Z}W_{n}^{[\alpha]}$ be the projection. For
any $u, v \in V$, $w \in W$ and $w^{'} \in W^{'}$, there exists a
multivalued analytic function of the form
\[
f(z_1, z_2) = \sum_{i,
j, k, l = 0}^N a_{ijkl}z_1^{m_i}z_2^{n_j}({\rm log}z_1)^k({\rm
log}z_2)^l(z_1 - z_2)^{-t}
\]
for $N \in \mathbb{N}$, $m_1, \dots,
m_N$, $n_1, \dots, n_N \in \mathbb{C}$ and $t \in \mathbb{Z}_{+}$,
such that the series
\[
\langle w^{'}, Y^{g; p}(u, z_1)Y^{g;p}(v,
z_2)w\rangle = \sum_{n \in \mathbb{C}}\langle w^{'}, Y^{g; p}(u,
z_1)\pi_nY^{g; p}(v, z_2)w\rangle,
\]
\[
\langle w^{'}, Y^{g; p}(v,
z_2)Y^{g;p}(u, z_1)w\rangle = \sum_{n \in \mathbb{C}}\langle w^{'},
Y^{g; p}(v, z_2)\pi_nY^{g; p}(u, z_1)w\rangle,
\]
\[
\langle w^{'}, Y^{g; p}(Y(u, z_1 - z_2)v,
z_2)w\rangle = \sum_{n \in \mathbb{C}}\langle w^{'}, Y^{g;
p}(\pi_nY(u, z_1 - z_2)v, z_2)w\rangle
\]
are absolutely convergent in
the regions $|z_1| > |z_2| > 0$, $|z_2| > |z_1| > 0$, $|z_2| > |z_1
- z_2| > 0$, respectively, and are convergent to the branch
\[
 \sum_{i, j, k, l = 0}^N
a_{ijkl}e^{m_il_p(z_1)}e^{n_jl_p(z_2)}l_p(z_1)^kl_p(z_2)^l(z_1 -
z_2)^{-t}
\]
of $f(z_1, z_2)$ when $\arg z_{1}$ and $\arg z_{2}$ are sufficiently close (more precisely, 
when $|\arg z_{1}-\arg z_{2}|<\frac{\pi}{2}$).

\item The {\it $\overline{\C}_{+}$- and $g$-grading conditions}: For $v\in V_{(m)}$
 and $w\in W_{p}=\coprod_{\alpha\in \C/\Z}W_{p}^{[\alpha]}$ where 
$m\in \Z$ and $p\in \overline{\C}_{+}$, write
$Y^{g}(v, x)w=\sum_{k=0}^{N}\sum_{n\in \C}Y^{g}_{n, k}(v)wx^{-n-1}(\log x)^{k}.$
Then $Y^{g}_{n, k}(v)w$ is $0$ when $m-n-1+p\not\in \overline{\C}_{+}$ (or $\Re(m-n-1+p)<0$),
is  in $W_{m-n-1+p}$ when $m-n-1+p\in \overline{\C}_{+}$ and for $r\in \R$,
$\coprod_{n\in \mathbb{I}}W_{n}$ is equal to the subspace of $W$ consisting 
of $w\in W$ such that when  $Y^{g}_{n, k}(v)w=0$ for 
$v\in V_{(m)}$, $n\in \C$, $k=0, \dots, N$, $m-n-1\not\in \overline{\C}_{+}$. 
For $\alpha \in \mathbb{C}/\mathbb{Z}$, $w \in
W^{[\alpha]}=\coprod_{n\in \overline{\C}_{+}}W_{n}^{[\alpha]}$, there exists $\Lambda \in \mathbb{Z}_{+}$
such that $ (g -
e^{2\pi\sqrt{-1}\alpha})^{\Lambda}w = 0$. Moreover, $gY^{g}(u, x)v=Y^{g}(gu, x)gv$.

\item The $L(-1)$-{\it derivative property}: For $v \in V$,
\[
\frac{d}{dx}Y^g(v, x) = Y^g(L(-1)v, x).
\]

\end{enumerate}
A {\it lower bounded generalized $g$-twisted $V$-module} or simply a {\it lower bounded $g$-twisted $V$-module}
is a $\overline{\C}_{+}$-graded weak $g$-twisted
$V$-module $W$ together with a decomposition of $W$ as a direct sum $W=\coprod_{n\in \C}W_{[n]}$
of generalized eigenspaces $W_{[n]}$ with eigenvalues $n\in \C$ for the operator 
$L^{g}(0)=\res_{x}xY^g(\omega, x) $ such that
for each $n\in \C$ and each $\alpha\in \C/\Z$, $W^{[\alpha]}_{[n + l]} =W_{[n+l]}\cap W^{[\alpha]}= 0$ for
sufficiently negative real number $l$. A lower bounded generalized $g$-twisted $V$-module
is said to be {\it strongly $\C/\Z$-graded} or {\it grading-restricted} if it is lower 
bounded and  for each $n \in
\mathbb{C}$, $\alpha \in \mathbb{C}/\mathbb{Z}$, 
$\dim W^{[\alpha]}_{[n]} =\dim W_{[n]}\cap W^{[\alpha]}< \infty$.}
\end{defn}

In the twisted representation theory of vertex operator algebras, we are mainly interested in strongly
$\mathbb{C}/\mathbb{Z}$-graded or grading-restricted 
generalized $g$-twisted modules. For simplicity, we shall call them 
simply $g$-twisted $V$-modules. In the case that $g$ does not act on $V$ semisimply, 
since the twisted vertex operators involve not only powers of the variable but also the 
logarithm of the variable, we shall also call such a $g$-twisted $V$-module a 
{\it logarithmic $g$-twisted $V$-module}. 
Lower bounded $g$-twisted $V$-modules also often appear 
in the formulations of certain conditions and assumptions. $\overline{\C}_{+}$-weak $g$-twisted
$V$-modules are not really needed in the twisted representation theory of vertex operator algebras.
They are not introduced in \cite{H} because of this reason. 
But in this paper, they are needed  in the proof of the theorem that $Z_{g}(V)$ and $A_{g}(V)$ are isomorphic.
Because of this reason, we prove all the results in this paper for $\overline{\C}_{+}$-weak $g$-twisted
$V$-modules and then derive the corresponding results for lower bounded $g$-twisted $V$-modules
and $g$-twisted $V$-modules (including  logarithmic $g$-twisted $V$-modules) as special cases.

Let $W=\coprod_{n\in \overline{\C}_{+},\alpha\in \C/\Z}W_{n}^{[\alpha]}$ be a 
$\overline{\C}_{+}$-graded weak $g$-twisted
$V$-module with the twisted vertex operator map $Y^{g}$.
For $\alpha\in \C/\Z$, let $W^{[\alpha]}=\coprod_{n\in \overline{\C}_{+}}W_{n}^{[\alpha]}$. 
When $W$ is a (lower bounded) $g$-twisted $V$-module,  we also have $W^{[\alpha]}=
\coprod_{n\in \C}W_{[n]}^{[\alpha]}$. 
As in the case of $V$, we call $a\in [0, 1)+\mathbb{I}$ a
{\it $g$-weight of $W$} if $W^{[a+\Z]} \neq 0$. 
Let $P(W)$ be the set of all $g$-weights of $W$. Then
$W=\coprod_{n\in \overline{\C}_{+}, a\in P(W)}W_{[n]}^{[a+\Z]}$. Let $\A_{W}: W\to W$ be the linear map defined 
by $\A_{W} w=aw$ for $w\in W^{[a+\Z]}$ and $a\in P(W)$. 
It is clear that for $a\in P(V)$ and $b\in P(W)$, 
either $a+b\in P(W)$ or $a+b-1\in P(W)$.

For $v\in V$, we have
$$Y^{g}(v, x)=\sum_{k=0}^{N}\sum_{n\in \C}Y^{g}_{n, k}(v)x^{-n-1}(\log x)^{k}.$$
Let
$$Y^{g}_{k}(v, x)=\sum_{n\in \C}Y^{g}_{n, k}(v)x^{-n-1}.$$
Then
$$Y^{g}(v, x)=\sum_{k=0}^{N}Y^{g}_{k}(v, x)(\log x)^{k}.$$
Denote the formal variable $\log x$ by $y$ and let
\begin{equation}\label{y-v-x-y}
Y^{g}(v, x, y)=\sum_{k=0}^{N}Y^{g}_{k}(v, x)y^{k}=\sum_{k=0}^{N}
\sum_{n\in \C}Y^{g}_{n, k}(v)x^{-n-1}y^{k}
\end{equation}
for $v\in V$.
Taking $u$ in the duality property to be $\one$, we see that $Y^{g}_{n, k}(v)w=0$ for $k=0,
\dots, N$ when 
$\Re(n)$ is sufficiently negative, that is, $Y^{g}_{k}(v, x)w$ for $k=0,
\dots, N$ are lower truncated. 
We say that $Y^{g}_{k}$, $k=0,
\dots, N$, are {\it lower truncated}.

Since the homogeneous subspaces of $V$ are finite dimensional and $g$ preserves the
homogeneous subspaces of $V$, by the multiplicative Jordan-Chevalley decomposition, there
exist a unique semisimple automorphism $\sigma$ of $V$ and a unique
locally unipotent operator $g_{u}$
such that $\sigma$ and $g_{u}$ commute with each other and $g=\sigma g_{u}$. Here
by locally unipotent operator, we mean that $g_{u}-1_{V}$ is locally nilpotent, that is,
for $v\in V$, $(g_{u}-1_{V})^{j}v=0$ when $j$ is sufficiently large.
Let
\begin{eqnarray}\label{def-N}
\mathcal{N}&=&\frac{\log g_{u}}{2\pi \sqrt{-1}}\nn
&=&\frac{\log (1_{V}+(g_{u}-1_{V}))}{2\pi \sqrt{-1}}\nn
&=&\sum_{j\in\Z_{+}}\frac{(-1)^{j+1}(g_{u}-1_{V})^{j}}{2\pi\sqrt{-1} j}.
\end{eqnarray}
Note that since $g_{u}$ is locally unipotent,  the right-hand side of (\ref{def-N})
is a finite sum when acting on $v\in V$ and is locally nilpotent. From the definition of $\mathcal{N}$, we have
\begin{eqnarray*}
g_{u}&=&e^{2\pi \sqrt{-1} \mathcal{N}}\nn
&=&\sum_{j\in \N}\frac{(2\pi \sqrt{-1}\mathcal{N})^{j}}{j!}.
\end{eqnarray*}
Thus we have $g=\sigma e^{2\pi \sqrt{-1} \mathcal{N}}$.

Let $(W, Y^{g})$ be a $\overline{\C}_{+}$-graded weak $g$-twisted
$V$-module. We now use the  multiplicative
Jordan-Chevalley decomposition
of $g$ above to derive some formulas.
For $n\in \Z$, consider the space $Y^{g}(V_{(n)}, x, y)$
of elements of $({\rm End}\;W)\{x\}[y]$ of the form
$Y^{g}(v, x, y)$ for $v\in V_{(n)}$. Then $Y^{g}(\cdot, x, y)$ is a surjective
homomorphism from $V_{(n)}$ to $Y^{g}(V_{(n)}, x, y)$.
In particular, $Y^{g}(V_{(n)}, x, y)$ is finite dimensional.
Since $g$ preserve weights, $g$ acts on $V_{(n)}$. Thus $g$ also acts on
$Y^{g}(V_{(n)}, x, y)$ by
$g\cdot Y^{g}(v, x, y)=Y^{g}(gv, x, y)$ for $v\in V_{(n)}$. Since on $V_{(n)}$,
$g=\sigma e^{2\pi \sqrt{-1}\mathcal{N}}$, we obtain
a decomposition of this action
of $g$ on $Y^{g}(V_{(n)}, x, y)$ which, for simplicity, we shall still write as
$g=\sigma e^{2\pi \sqrt{-1}\mathcal{N}}$
where $\sigma$ and $\mathcal{N}$ are the corresponding operators on $Y^{g}(V_{(n)}, x, y)$.

On the other hand, from the equivariance property,
\begin{equation}\label{monodromy}
Y^{g}(gv, x, y) = e^{-2\pi \sqrt{-1} x\partial/\partial x}e^{-2 \pi \sqrt{-1} \partial/\partial y}Y^{g}(v, x, y)
\end{equation}
for $v\in V$.
By (\ref{monodromy}), the action of $g$ on $Y^{g}(V_{(n)}, x, y)$ is equal to
$e^{-2\pi\sqrt{-1}x\partial/\partial x}e^{-2 \pi \sqrt{-1} \partial/\partial y}$. If $Y^{g}(V_{(n)}, x, y)$
is invariant under $e^{-2\pi \sqrt{-1} x\partial/\partial x}$ and
$e^{-2 \pi\sqrt{-1} \partial/\partial y}$, then $e^{-2\pi \sqrt{-1} x\partial/\partial x}e^{-2 \pi\sqrt{-1} \partial/\partial y}$ is
a Jordan-Chevalley decomposition of $g$ on $Y^{g}(V_{(n)}, x, y)$.
By the uniqueness of  the Jordan-Chevalley decomposition, we obtain
$\sigma=e^{-2\pi\sqrt{-1}x\partial/\partial x}$ and $e^{2\pi \sqrt{-1} \mathcal{N}}=e^{-2 \pi \sqrt{-1} \partial/\partial y}$.
Since
\begin{eqnarray*}
\partial/\partial y&=&\frac{\log (1_{Y^{g}(V_{(n)}, x, y)}
+(e^{-2 \pi\sqrt{-1} \partial/\partial y}-1_{Y^{g}(V_{(n)}, x, y)})}{-2\pi \sqrt{-1}}\nn
&=&\sum_{j\in\Z_{+}}\frac{(-1)^{j+1}(e^{-2 \pi\sqrt{-1} \partial/\partial y}
-1_{Y^{g}(V_{(n)}, x, y)})^{j}}{-2\pi \sqrt{-1} j},
\end{eqnarray*}
by the definition of $\mathcal{N}$, we obtain
$\mathcal{N}=-\partial/\partial y$ on $Y^{g}(V_{(n)}, x, y)$.
Since we do not assume that $Y^{g}(V_{(n)}, x, y)$
is invariant under $e^{-2\pi \sqrt{-1} x\partial/\partial x}$ and
$e^{-2 \pi\sqrt{-1} \partial/\partial y}$, we have to prove 
$\sigma=e^{-2\pi\sqrt{-1} x\partial/\partial x}$ and $e^{2\pi \sqrt{-1} \mathcal{N}}=e^{-2 \pi\sqrt{-1} \partial/\partial y}$
directly, without using the uniqueness of  the Jordan-Chevalley decomposition for an operator 
on a finite-dimensional space.

Let $v\in V^{[\alpha]}_{(n)}$, $\alpha\in \C/\Z$. Since the action of $g$ on $Y^{g}(V_{(n)}, x, y)$ is equal to
$e^{-2\pi \sqrt{-1} x\partial/\partial x}e^{-2 \pi \sqrt{-1} \partial/\partial y}$ and $Y^{g}(v, x, y)$ is a generalized 
eigenvector of $g$ with eigenvalue $e^{2\pi \sqrt{-1}\alpha}$, there exists $M\in \Z_{+}$ such that 
\begin{equation}\label{j-c-decomp-1}
(e^{-2\pi \sqrt{-1}x\partial/\partial x}e^{-2 \pi\sqrt{-1} \partial/\partial y}-e^{2\pi \sqrt{-1}\alpha})^{M}Y^{g}(v, x, y)=0.
\end{equation}
From (\ref{y-v-x-y}) and (\ref{j-c-decomp-1}), we obtain
\begin{eqnarray}\label{j-c-decomp-2}
\lefteqn{\sum_{k=0}^{N}
\sum_{m\in \C}Y^{g}_{m, k}(v)(e^{-2\pi\sqrt{-1}(-m-1)}e^{-2 \pi\sqrt{-1}\partial/\partial y}-e^{2\pi \sqrt{-1}\alpha})^{M}x^{-m-1}y^{k}}\nn
&&=\sum_{k=0}^{N}
\sum_{m\in \C}Y^{g}_{m, k}(v)(e^{-2\pi\sqrt{-1} x\partial/\partial x}e^{-2 \pi\sqrt{-1}\partial/\partial y}-e^{2\pi\sqrt{-1}\alpha})^{M}x^{-m-1}y^{k}\nn
&&=0.
\end{eqnarray}
From (\ref{j-c-decomp-2}), we see that the coefficients of the left-hand side of (\ref{j-c-decomp-2})
as a formal series in powers of $x$ must be $0$, that is,
\begin{equation}\label{j-c-decomp-3}
\sum_{k=0}^{N}
Y^{g}_{m, k}(v)(e^{-2\pi \sqrt{-1} (-m-1)}e^{-2 \pi\sqrt{-1} \partial/\partial y}-e^{2\pi \sqrt{-1}\alpha})^{M}y^{k}=0
\end{equation}
for $m\in \C$. Rewriting (\ref{j-c-decomp-3}), we obtain
$$
(e^{-2 \pi \sqrt{-1}\partial/\partial y}-e^{2\pi \sqrt{-1}(\alpha-m-1)})^{M}\sum_{k=0}^{N}
Y^{g}_{m, k}(v)y^{k}=0
$$
for $m\in \C$. Let $m\in \C$ such that  $e^{2\pi\sqrt{-1}(\alpha-m-1)}\ne 1$. If $Y^{g}_{m, k}(v)\ne 0$
for some $k$, 
then $\sum_{k=0}^{N}
Y^{g}_{m, k}(v)y^{k}$
is a nonzero polynomial in $y$. Since $e^{2\pi\sqrt{-1}i(\alpha-m-1)}\ne 1$,
$$(e^{-2 \pi \sqrt{-1} \partial/\partial y}-e^{2\pi \sqrt{-1}(\alpha-m-1)})^{M}\sum_{k=0}^{N}
Y^{g}_{m, k}(v)y^{k}$$
cannot be equal to $0$. Contradiction with (\ref{j-c-decomp-2}).
Thus for such $m$, $Y^{g}_{m, k}(v)= 0$ for $k=0, \dots, N$ and we have 
$$
Y^{g}(v, x, y)=\sum_{k=0}^{N}
\sum_{m\in \alpha}Y^{g}_{m, k}(v)x^{-m-1}y^{k}.
$$
In particular, 
\begin{eqnarray*}
e^{-2\pi \sqrt{-1} x\partial/\partial x}Y^{g}(v, x, y)&=&e^{2\pi \sqrt{-1} \alpha}Y^{g}(v, x, y)\nn
&=&\sigma Y^{g}(v, x, y).
\end{eqnarray*}
Since $v$ is an arbitrary element of $V_{(n)}$,  $e^{-2\pi \sqrt{-1} x\partial/\partial x}=\sigma$
on $Y^{g}(V_{(n)}, x, y)$ and thus $e^{2\pi\sqrt{-1} \mathcal{N}}=e^{-2 \pi \sqrt{-1} \partial/\partial y}$.

Since $n$ is arbitrary, we have proved the following formulas:

\begin{lemma}
For $v\in V$,
\begin{eqnarray}
Y^{g}(\sigma v, x, y) &=& e^{-2\pi ix\partial/\partial x}Y^{g}(v, x, y),\label{sigma=xdx}\\
Y^{g}(\mathcal{N}v, x, y) &=&-\partial/\partial yY^{g}(v, x, y),\label{n=dy}\\
Y^{g}(e^{2\pi i\mathcal{N}}v, x, y)& =&e^{-2\pi i \partial/\partial y}Y^{g}(v, x, y).\label{en=edy}\epfe
\end{eqnarray}
\end{lemma}

From (\ref{sigma=xdx}), (\ref{n=dy}) and (\ref{en=edy}), for $v\in V^{[a+\Z]}$, $a\in P(V)$ and
$w\in W$, $Y^{g}(v, x)w\in x^{-a}W((x))[\log x]$.
Then for $v \in V^{[\alpha]}$, $w \in W$,
\begin{eqnarray*}
Y^g_k(v, x)w &=& \sum_{n \in a + \Z}Y^{g}_{n,k}(v)wx^{-n-1},\\
Y^g(v, x)w &=& \sum_{k = 0}^{N_g(v)}Y^g_k(v, x)w(\log x)^k,
\end{eqnarray*}
where $N_g(v)\in \N$ depends only on $\mathcal{N}$ and $v$.
For convenience, we shall use $u(n)$ to denote the operator $Y^{g}_{n,0}(u)$.

From (\ref{n=dy}), we have the following relation between $Y^{g}$ and  $Y^{g}_{0}$:
\begin{lemma}\label{lemma-yg-y0}
For $v\in V$,
\begin{equation}\label{yg-y0}
Y^g(v, x) = Y^{g}_{0}(x^{-\mathcal{N}}v, x),
\end{equation}
where $x^{-\mathcal{N}} = e^{-\mathcal{N}\log x}$.
\end{lemma}
\pf
From (\ref{n=dy}) and the formal Taylor's theorem, we obtain
\begin{eqnarray}\label{yg-y0-1}
Y^g(x^{\mathcal{N}}v, x, y) &=&e^{-(\log x)\partial/\partial y}Y^{g}(v, x, y)\nn
&=&\sum_{k=0}^{N_g(v)}e^{-(\log x)\partial/\partial y}Y^{g}_{k}(v, x)y^{k}\nn
&=&\sum_{k=0}^{N_g(v)}Y^{g}_{k}(v, x)(y-\log x)^{k}.
\end{eqnarray}
Substituting $\log x$ for $y$ in (\ref{yg-y0-1}) and then substituting
$x^{-\mathcal{N}}v$ for $v$, we obtain (\ref{yg-y0}).
\epfv

We also need the following lemma:
\begin{lemma}\label{sigma-N-Y}
For $u\in V$, we have
\begin{eqnarray}
\sigma Y^{g}(u, x)&=&Y^{g}(\sigma u, x)\sigma,\label{sigma-Y}\\
{[\mathcal{N}, Y^{g}(u, x)]}&=&Y^{g}(\mathcal{N}u, x),\label{N-Y}\\
e^{y\mathcal{N}}Y^{g}(u, x)&=&Y^{g}(e^{y\mathcal{N}}u, x)e^{y\mathcal{N}}.\label{e-N-Y}
\end{eqnarray}
\end{lemma}
\pf
Let $u\in V^{[\alpha]}$ and $v\in V^{[\beta]}$. Then $u$ and $v$ are generalized eigenvectors
for $g$ with eigenvalues $e^{2\pi\sqrt{-1}\alpha}$ and $e^{2\pi \sqrt{-1}\beta}$.
In particular, there exist $M, N\in \Z_{+}$ such that
$(g-e^{2\pi\sqrt{-1}\alpha})^{M}u=(g-e^{2\pi\sqrt{-1}\beta})^{N}v=0$.

Using the fact that $g$ is an automorphism of $V$, we have
$$(g-e^{2\pi \sqrt{-1} (\alpha+\beta)})Y^{g}(u, x)v
=Y^{g}((g-e^{2\pi \sqrt{-1}\alpha})u, x)gv+Y^{g}(e^{2\pi \sqrt{-1} \alpha}u, x)(g-e^{2\pi \sqrt{-1} \beta})v.$$
Then
\begin{eqnarray*}
\lefteqn{(g-e^{2\pi \sqrt{-1} (\alpha+\beta)})^{M+N}Y^{g}(u, x)v}\nn
&&=\sum_{j=0}^{M+N}{M+N\choose j}Y^{g}(e^{2\pi \sqrt{-1} (M+N-j)\alpha}(g-e^{2\pi \sqrt{-1} \alpha})^{j}u, x)
g^{j}(g-e^{2\pi \sqrt{-1} \beta})^{M+N-j}v\nn
&&=0.
\end{eqnarray*}
Thus $Y^{g}(u, x)v\in V^{[\alpha+\beta]}[[x, x^{-1}]]$. This proves (\ref{sigma-Y}).

Using (\ref{sigma-Y}) together with the fact that $g$ is an automorphism of $V$,
we obtain
$$e^{2\pi\sqrt{-1}\mathcal{N}}Y^{g}(u, x)=Y^{g}(e^{2\pi \sqrt{-1}\mathcal{N}}u, x)e^{2\pi \sqrt{-1}\mathcal{N}}$$
or
\begin{eqnarray*}
Y^{g}(e^{2\pi\sqrt{-1}\mathcal{N}}u, x)&=&e^{2\pi \sqrt{-1}\mathcal{N}}Y^{g}(u, x)e^{-2\pi\sqrt{-1}\mathcal{N}}\nn
&=&{\rm Ad}(e^{2\pi\sqrt{-1}\mathcal{N}})(Y^{g}(u, x)).
\end{eqnarray*}
Thus by the definition of $\mathcal{N}$,
\begin{eqnarray}\label{N-Y-1}
Y^{g}(\mathcal{N}u, x)&=&\frac{1}{2\pi \sqrt{-1}}Y^{g}(\log  (e^{2\pi \sqrt{-1}\mathcal{N}})u, x)\nn
&=&\frac{1}{2\pi \sqrt{-1}}(\log  ({\rm Ad}(e^{2\pi \sqrt{-1}\mathcal{N}}))))(Y^{g}(u, x))\nn
&=&\frac{1}{2\pi \sqrt{-1}}{\rm ad}(2\pi \sqrt{-1}\mathcal{N})(Y^{g}(u, x))\nn\nn
&=&[\mathcal{N}, Y^{g}(u, x)],
\end{eqnarray}
proving (\ref{N-Y}). The formula (\ref{e-N-Y})  follows immediately from (\ref{N-Y}).
\epfv

As consequences of Lemmas \ref{lemma-yg-y0} and  \ref{sigma-N-Y}, we have the following properties:
\begin{prop}
For $n \in \Z$ and $v\in V$,
\begin{eqnarray}
[\mathcal{N}, L(n)] &=&0,\label{bracket-N-L}\\
x^{\mathcal{N}}L(n) &= &L(n)x^{\mathcal{N}},\label{conj-N-L}\\
Y_0^g(L(-1)u, x) &=& \frac{d}{dx}Y_0^g(u, x)-x^{-1}Y_0^g(\mathcal{N}u, x).\label{L(-1)derivative}
\end{eqnarray}
\end{prop}
\pf
Since $g \omega=\omega$, we have $\mathcal{N}\omega=\omega$. Then (\ref{bracket-N-L})
follows from (\ref{N-Y}) and (\ref{conj-N-L}) is an immediate consequence of  (\ref{bracket-N-L}).

The formula (\ref{L(-1)derivative}) is obtained by taking derivative with respect to $x$ on
both sides of (\ref{yg-y0}), using
the $L(-1)$-derivative property for $Y^{g}$ and  then take the $0$-th power of $\log x$
on both sides of the resulting equality.
\epfv

Using the the formulas obtained above, we prove a duality property for 
the map $Y_{0}^{g}$. For $u\in V$, in the region $|z_{2}|>|z_{1}-z_{2}|>0$, we use 
\begin{equation}\label{correction}
\left(1 + \frac{z_{1}-z_{2}}{z_2}\right)^{\mathcal{N}}u
\end{equation}
to denote the power series in the variable $\frac{z_{1}-z_{2}}{z_{2}}$ with coefficients in $V$ 
obtained by using the binomial expansion. This series can also be obtained by 
writing (\ref{correction}) as 
\begin{equation}\label{correction2}
\exp\left(\log \left(1 + \frac{z_{1}-z_{2}}{z_2}\right)\mathcal{N}\right)u,
\end{equation}
expanding (\ref{correction2}) as a power series in 
\begin{equation}\label{log}
\log \left(1 + \frac{z_{1}-z_{2}}{z_2}\right)
\end{equation}
and then expanding powers of  (\ref{log}) as powers series in 
$\frac{z_{1}-z_{2}}{z_{2}}$ in the region $|z_{2}|>|z_{1}-z_{2}|>0$. 
Since $\mathcal{N}$ is locally nilpotent, this series is in fact obtained from a polynomial 
in (\ref{log}) with coefficients in $V$ by expanding powers of  (\ref{log}) as powers series in 
$\frac{z_{1}-z_{2}}{z_{2}}$ in $|z_{2}|>|z_{1}-z_{2}|>0$. 
In other words,
(\ref{correction}) is a finite sum of elements of $V$ multiplied by powers of (\ref{log}), which
as power series in $\frac{z_{1}-z_{2}}{z_2}$ are absolutely convergent
in the region $|z_{2}|>|z_{1}-z_{2}|>0$.

\begin{thm}\label{duality-y-0}
 Let $W$ be a $\overline{\C}_{+}$-graded weak $g$-twisted
$V$-module. Then for $u \in V^{[a+\Z]}$,  $v\in V^{[b+\Z]}$,  $a, b\in P(V)$,  
\begin{eqnarray}
&\langle w', Y^{g}_{0}(u, z_1)Y^{g}_{0}(v, z_2)w\rangle,&\label{jacobi1}\\
&\langle w', Y^{g}_{0}(v, z_2)Y^{g}_{0}(u, z_1)w\rangle,&\label{jacobi2}\\
&{\ds \left\langle w', Y^{g}_{0}\left(Y\left(\left(1 + \frac{z_{1}-z_{2}}{z_2}\right)^{\mathcal{N}}u,
z_{1}-z_{2}\right)v, z_2\right)w\right\rangle}&
\label{jacobi3}
\end{eqnarray}
are absolutely convergent in the region $|z_{1}|>|z_{2}|>0$, $|z_{2}|>|z_{1}|>0$,
$|z_{2}|>|z_{1}-z_{2}|>0$, respectively, to
a function of $z_{1}$ and $z_{2}$ of the form
\begin{equation}\label{g-rat-fn}
e^{-a\log z_{1}}e^{-b \log z_{2}}\frac{g(z_1, z_2)}{z_{1}^{r} z_{2}^{s}(z_{1}-z_{2})^{t}},
\end{equation}
where $g(z_1, z_2)$ is a polynomial in $z_{1}$ and $z_{2}$, $r, s, t \in \Z$.
\end{thm}
\pf
Since $u \in V^{[a+\Z]}$, $v\in V^{[b+\Z]}$, $a, b\in P(V)$,   by the duality property,
\begin{eqnarray}
&\langle w', Y^{g}(u, z_1)Y^{g}(v, z_2)w\rangle,&\label{prod12}\\
&\langle w', Y^{g}(v, x_2)Y^{g}(u, x_1)w\rangle,&\label{prod21}\\
&{\ds \langle w', Y^{g}(Y(u, z_{1}-z_{2})v, z_2)w\rangle}&\label{iter}
\end{eqnarray}
are absolutely convergent in the region $|z_{1}|>|z_{2}|>0$, $|z_{2}|>|z_{1}|>0$,
$|z_{2}|>|z_{1}-z_{2}|>0$, respectively, and are convegent when $|\arg z_{1}-\arg z_{2}|<\frac{\pi}{2}$ to
a function of $z_{1}$ and $z_{2}$ of the form
\begin{eqnarray*}
\lefteqn{ \sum_{i, j, k, l = 0}^N
a_{ijkl}e^{(-a-r_{i})\log z_1}e^{(-b-s_{j})\log z_2}(\log z_1)^k (\log z_2)^l(z_1 -
z_2)^{-t}}\nn
&&=e^{-a\log z_{1}}e^{-b \log z_{2}}\sum_{k, l=0}^{N}\frac{g_{kl}(z_{1}, z_{2})}{z_{1}^{r_{0}}z_{2}^{s_{0}}(z_1 -
z_2)^{t}}(\log z_1)^k (\log z_2)^l
\end{eqnarray*}
where $r_{i}, s_{j}, t \in \Z$ satisfying $r_{i}<r_{i-1}$ and $s_{j}<s_{j-1}$ and $g_{kl}(z_{1}, z_{2})$
are polynomials in $z_{1}$ and $z_{2}$. In particular, in the region $|z_{1}|>|z_{2}|>0$,
\begin{eqnarray*}
\lefteqn{\sum_{k, l=0}^{N}\langle w', Y^{g}_{k}(u, z_1)Y^{g}_{l}(v, z_2)w\rangle(\log z_1)^k (\log z_2)^l}\nn
&&=\langle w', Y^{g}(u, z_1)Y^{g}(v, z_2)w\rangle\nn
&&=e^{-a\log z_{1}}e^{-b \log z_{2}}\sum_{k, l=0}^{N}\frac{g_{kl}(z_{1}, z_{2})}{z_{1}^{r_{0}}z_{2}^{s_{0}}(z_1 -
z_2)^{t}}(\log z_1)^k (\log z_2)^l.
\end{eqnarray*}
By Proposition 7.8 in \cite{HLZ5}, we obtain
$$\langle w', Y^{g}_{k}(u, z_1)Y^{g}_{l}(v, z_2)w\rangle=
e^{-a\log z_{1}}e^{-b \log z_{2}}\frac{g_{00}(z_{1}, z_{2})}{z_{1}^{r_{0}}z_{2}^{s_{0}}(z_1 -
z_2)^{t}}.$$
Taking $g(z_{1}, z_{2})$ to be $g_{00}(z_{1}, z_{2})$ and $r, s$ to be $r_{0}, s_{0}$,
we see that (\ref{jacobi1}) indeed converges absolutely to a function of the form (\ref{g-rat-fn}).

Similarly we can prove that (\ref{jacobi2}) converges absolutely to a function of the form (\ref{g-rat-fn}).
Since (\ref{prod12}) and (\ref{prod21}) converge to the same function in $|z_{1}|>|z_{2}|>0$ and $|z_{2}|>|z_{1}|>0$,
respectively, from the argument above, the fact that the powers of $z_{1}$ and $z_{2}$ in these expressions are
congruent to  $a$ and $b$, respectively, and Proposition 7.8 in \cite{HLZ5}, we see that
(\ref{jacobi1}) and (\ref{jacobi2}) converge absolutely to the same function
of the form (\ref{g-rat-fn})  in $|z_{1}|>|z_{2}|>0$ and $|z_{2}|>|z_{1}|>0$, respectively.

The same argument  above shows that
$$\langle w', Y^{g}_{0}(Y(u,
z_{1}-z_{2})v, z_2)w\rangle$$
converges absolutely in the region $|z_{2}|>|z_{1}-z_{2}|>0$.
As is discussed above, (\ref{correction}) is a finite sum of elements of $V$ multiplied by powers of (\ref{log}),
which as power series in $\frac{z_{1}-z_{2}}{z_2}$ are absolutely convergent
in the region $|z_{2}|>|z_{1}-z_{2}|>0$,
(\ref{jacobi3}) is absolutely convergent in the same region.
We shall use the associativity for $Y^{g}$ and (\ref{yg-y0}) to show that (\ref{jacobi3}) converges to
the same function that (\ref{jacobi1}) converges to.

By  (\ref{yg-y0}), (\ref{e-N-Y}) and the associativity for $Y^{g}$, in the region $|z_{1}|>|z_{2}|>|z_{1}-z_{2}|>0$,
we have
\begin{eqnarray*}
\lefteqn{{\ds \left\langle w', Y^{g}_{0}\left(Y\left(\left(1 + \frac{z_{1}-z_{2}}{z_2}\right)^{\mathcal{N}}u,
z_{1}-z_{2}\right)v, z_2\right)w\right\rangle}}\nn
&&={\ds \left\langle w', Y^{g}\left(z_{2}^{\mathcal{N}}Y\left(\left(\frac{z_{1}}{z_2}\right)^{\mathcal{N}}u,
z_{1}-z_{2}\right)v, z_2\right)w\right\rangle}\nn
&&={\ds \langle w', Y^{g}(Y(z_{1}^{\mathcal{N}}u,
z_{1}-z_{2})z_{2}^{\mathcal{N}}v, z_2)w\rangle}\nn
&&=\langle w', Y^{g}((z_{1})^{\mathcal{N}}u,
z_{1})Y^{g}(z_{2}^{\mathcal{N}}v, z_2)w\rangle\nn
&&=\langle w', Y^{g}_{0}(u,
z_{1})Y^{g}_{0}(v, z_2)w\rangle.
\end{eqnarray*}
\epfv

In \cite{B}, Bakalov derived a Jacobi identity for the map $Y_{0}^{g}$ from an associator formula
for twisted $V$-modules. From Theorem \ref{duality-y-0}, we obtain 
the Jacobi identity immediately.

\begin{thm}[Bakalov \cite{B}]
Let $W$ be a $\overline{\C}_{+}$-graded weak $g$-twisted
$V$-module. Then  we have the twisted Jacobi identity
\begin{eqnarray}\label{Jacobi}
x_0^{-1}\delta\left(\frac{x_1 - x_2}{x_0}\right)Y^{g}_{0}(u, x_1)Y^{g}_{0}(v, x_2)
- x_0^{-1}\delta\left(\frac{- x_2 + x_1}{x_0}\right)Y^{g}_{0}(v, x_2)Y^{g}_{0}(u, x_1)\nn
= x_1^{-1}\delta\left(\frac{x_2+x_0}{x_1}\right)\left(\frac{x_2+x_0}{x_1}\right)^{a}
Y^{g}_{0}\left(Y\left(\left(1 + \frac{x_0}{x_2}\right)^{\mathcal{N}}u, x_0\right)v, x_2\right)
\end{eqnarray}
or equivalently, the component form of the twisted Jacobi identity 
\begin{eqnarray}\label{jacobi-comp}
&{\ds \sum_{k\in \N}(-1)^{k}\binom{l}{k}u(l-k+m)v(k+n)
-\sum_{k\in \N}(-1)^{k-l}\binom{l}{k}v(l-k+n)u(k+m)}&\nn
&{\ds = \sum_{j \in \N} \left(\left(\binom{l+\mathcal{N}}{j}u\right)_{j}
v\right)(m+n-j)}&
\end{eqnarray}
for $u\in V^{[a+\Z]}$, $v\in V^{[b+\Z]}$, $m\in a+\Z$, $n\in b+\Z$ and $l\in \Z$,
where $u(p)=Y_{p, 0}^{g}(u)$ and $v(q)=Y_{q, 0}^{g}(v)$ for $p\in a+\Z$ and $q\in b+\Z$. 
\epf
\end{thm}


As a consequence of Jacobi identity, we have the commutator formulas:
\begin{cor}
For $u \in V^{[a+\Z]}$, $v\in V^{[b+\Z]}$, $m\in a+\Z$ and $n\in b+\Z$, we have
\begin{eqnarray}\label{commu}
{[Y^{g}_{0}(u, x_1), Y^{g}_{0}(v, x_2)] }
&=& \res_{x_0}x_1^{-1}\delta\left(\frac{x_2+x_0}{x_1}\right)
\left(\frac{x_2+x_0}{x_1}\right)^{a}\cdot\nn
&&\quad\quad\quad\cdot  Y^{g}_{0}\left(Y\left(\left(1+\frac{x_0}{x_2}\right)
^{\mathcal{N}}u, x_0\right)v, x_2\right),\label{comm3}\\
{[u(m), Y^{g}_{0}(v, x_2)] }&=& \sum_{j =0}^{\infty}x_2^{m-j}
Y^{g}_{0}\left(\left(\binom{m+\mathcal{N}}{j}u\right)(j)v, x_2\right),
\label{comm1}\\
{[u(m), v(n)] }&=& \sum_{j = 0}^{\infty} \left(\left(\binom{m+\mathcal{N}}{j}u\right)_{j}
v\right)(m+n-j).\label{comm2}
\end{eqnarray}
\end{cor}

We also have the associator formula and weak associativity:
\begin{cor}
For $u \in V^{[a+\Z]}$, $v\in V$ and $w \in W$, 
\begin{eqnarray}\label{associator}
&{\displaystyle \res_{x_{1}}x_{1}^{a}x_0^{-1}\delta\left(\frac{x_1 - x_2}{x_0}\right)Y^{g}_{0}(u, x_1)Y^{g}_{0}(v, x_2)
- \res_{x_{1}}x_{1}^{a} x_0^{-1}\delta\left(\frac{- x_2 + x_1}{x_0}\right)Y^{g}_{0}(v, x_2)Y^{g}_{0}(u, x_1)}&\nn
&{\displaystyle =(x_2+x_0)^{a}
Y^{g}_{0}\left(Y\left(\left(1 + \frac{x_0}{x_2}\right)^{\mathcal{N}}u, x_0\right)v, x_2\right).}
\end{eqnarray}
For $u \in V^{[a+\Z]}$, $v\in V$ and $w \in W$, let $l \in a + \Z$
such that $u(n)w = 0$ for $n \geq l$. Then
\begin{equation}\label{asso}
(x_0+x_2)^l Y^{g}_{0}(u, x_0+x_2)Y^{g}_{0}(v, x_2)w
= (x_2+x_0)^l Y^{g}_{0}\left(Y\left(\left(1+\frac{x_0}{x_2}\right)^{\mathcal{N}}u,
x_0\right)v, x_2\right)w.
\end{equation}
\end{cor}

Finally we have the following equivalence of the main properties of $\overline{\C}_{+}$-graded weak 
$g$-twisted modules:

\begin{thm}\label{jacobi-equiv}
The following properties for a $\overline{\C}_{+}$-graded weak $g$-twisted
$V$-module are equivalent:
\begin{enumerate}

\item[(1)] The duality property for $Y^{g}$ in Definition \ref{d}.

\item[(2)] The duality property for $Y^{g}_{0}$  in Theorem \ref{duality-y-0}.

\item[(3)] The property that $Y^{g}_{0}$ is lower truncated
and  the twisted Jacobi identity (\ref{Jacobi}).

\item[(4)] The property that $Y^{g}_{0}$ is lower truncated, the commutator formula (\ref{comm3}) and the weak associativity (\ref{asso}).
\end{enumerate}
\end{thm}
\pf We have shown that $(1) \Rightarrow (2) \Rightarrow (3) \Rightarrow (4)$. 
The same method as that for ordinary modules shows  that $(4) \Rightarrow (2)$.
Using (\ref{yg-y0}), we see that $(2)\Rightarrow (1)$.
\epfv

\setcounter{equation}{0}
\section{The $g$-twisted universal enveloping algebra and $g$-twisted zero-mode algebra}

In this section,  we construct two associative algebras associated with $V$ and the automorphism $g$,
one is called  the ``$g$-twisted universal algebra of $V$" and the other is called the ``$g$-twisted zero-mode algebra
associated to $V$".

Let 
$$L(V, g) = \coprod_{a \in P(V)} V^{[a+\Z]}\otimes t^{a}\C[t, t^{-1}].$$
For $u \in V^{[a+\Z]}$ and $m \in a + \Z$, we use $u(m)$ to denote the 
element $u\otimes t^m \in L(V, g)$.

Let $T(L(V,g))$ be the tensor algebra of $L(V, g)$. Then $T(L(V,g))$ is spanned by elements of the form
$$
u_1(i_1) \otimes \cdots \otimes u_n(i_n)$$
for $u_s \in V^{[a_{s}+\Z]},\; i_s \in a_{s} + \Z$.
For simplicity, we omit the tensor symbols for the elements. Then the multiplication is given by
$$
(u_1(i_1)\cdots u_m(i_m)) (v_1(j_1)\cdots v_n(j_n)) = u_1(i_1)\cdots u_m(i_m)v_1(j_1)\cdots v_n(j_n).
$$
Define a grading called {\it weight} on $T(L(V,g))$ as follows: For $\lambda\in \C$,
$$\wt \lambda=0$$
and 
for $\lambda\in \C$, homogeneous $u_s \in V^{[a_{s}+\Z]}$ and $i_{s}\in a_{s}+\Z$,
$$
\wt (\lambda u_1(i_1)\cdots u_n(i_n)) = \sum_{s=1}^n(\wt u_s- i_s - 1).
$$
Let $\tilde{P}(V)$ be the subset of $[0, 1)+\mathbb{I}$ such that $\tilde{P}(V)+\Z$ is the subgroup of $\C/\Z$ 
generated by $P(V)+\Z$. Then the weight of a homogeneous element of $T(L(V,g))$ is in $\tilde{P}(V)+\Z$. 
Let $T(L(V,g))_k$ denote the subspace of $T(L(V,g))$  consisting of the elements of weight $k$. Then we have
$$
T(L(V,g)) = \coprod_{k \in \tilde{P}(V)+\Z}T(L(V,g))_k
$$
and
$$
T(L(V,g))_k T(L(V,g))_l \subset T(L(V,g))_{k+l}.
$$

For $k\in \tilde{P}(V)+\Z$ and $n\in \Z$, let
$$
T(L(V,g))_{k}^{n}  = \coprod_{a\in \tilde{P}(V), i\in n+\N}T(L(V,g))_{k+a+i}  T(L(V,g))_{-a-i}.
$$
Then we have
$$T(L(V,g))_k^{n+1} \subset T(L(V,g))^{n}_k$$
and also
\begin{eqnarray*}
\bigcap_{n \in \Z}T(L(V,g))_k^{n}&=& 0, \\
\bigcup_{n \in \Z}T(L(V,g))_k^{n}& =& T(L(V,g))_k.
\end{eqnarray*}
Hence, $\{T(L(V,g))_k^n\;|\; n\in \Z\}$ is a fundamental neighborhood system of $T(L(V,g))_k$.
Denote by $\tilde{T}(L(V, g))_k$
its completion and let
\[
\tilde{T}(L(V, g)) = \coprod_{k \in \tilde{P}(V)+\Z}\tilde{T}(L(V, g))_k.
\]
Then $\tilde{T}(L(V, g))$ is a complete topological vector space. 

\begin{prop}
The multiplication for $T(L(V,g))$ is continuous under the topology given by the fundamental 
neighborhood system above. In particular, its topological completion $\tilde{T}(L(V, g))$ 
is a complete topological ring. 
\end{prop}
\pf
It suffices to prove that  given $m\in \tilde{P}(V)+ \Z$ and $n\in \Z$
and given $k, l\in \tilde{P}(V)+ \Z$  satisfying $k+l=m$,
there exist $n_{1}, n_{2}\in \Z$
 such that 
\begin{equation}\label{cont-1}
T(L(V, g))_{k}^{n_{1}}T(L(V, g))_{l}+T(L(V, g))_{k}T(L(V, g))_{l}^{n_{2}}\subset T(L(V, g))_{m}^{n}.
\end{equation}

Let $l=c+l'$ where $c\in \tilde{P}(V)$ and $l'\in \Z$. 
We take $n_{1}=n+l'+1$  and $n_{2}=n$. Elements of $T(L(V, g))_{k}^{n+l'}$
are finite sums of elements of the form $\varphi \psi$ where $\varphi\in T(L(V, g))_{k+a+i}$ 
and $\psi\in T(L(V, g))_{-a-i}$ for $a\in \tilde{P}(V)$, $i\in n+l'+1+\N$. Let $\eta\in T(L(V, g))_{l}$.
When $a-c\in \tilde{P}(V)$, $\varphi \in T(L(V, g))_{(k+l)+(a-c)+(i-l')}$ and 
$\psi\eta\in T(L(V, g))_{-(a-c)-(i-l')}$. Since $a-c\in \tilde{P}(V)$ and $i-l'\in n+1+\N\subset n+\N$, 
\begin{eqnarray*}
\varphi \psi \eta&\in& T(L(V, g))_{(k+l)+(a-c)+(i-l')}T(L(V, g))_{-(a-c)-(i-l')}\nn
&\subset &T(L(V, g))_{k+l}^{n}\nn
&=&T(L(V, g))_{m}^{n}.
\end{eqnarray*}
When $a-c\not\in \tilde{P}(V)$, we must have $a-c+1\in \tilde{P}(V)$.  Then 
$\varphi \in T(L(V, g))_{(k+l)+(a-c+1)+(i-l'-1)}$ and 
$\psi\eta\in T(L(V, g))_{-(a-c+1)-(i-l'-1)}$. Since $a-c+1\in \tilde{P}(V)$ and $i-l'-1\in n+\N$, 
\begin{eqnarray*}
\varphi \psi \eta& \in& T(L(V, g))_{(k+l)+(a-c+1)+(i-l'-1)}T(L(V, g))_{-(a-c+1)-(i-l'-1)}\nn
&\subset &T(L(V, g))_{k+l}^{n}\nn
&=&T(L(V, g))_{m}^{n}.
\end{eqnarray*}
Since the addition is continuous with respect to the topology and elements of $T(L(V, g))_{k}^{n+l'}$
are finite sums of elements of the form $\varphi \psi$ above, we have shown that 
\begin{equation}\label{cont-2}
T(L(V, g))_{k}^{n_{1}}T(L(V, g))_{l}\subset T(L(V, g))_{m}^{n}.
\end{equation}

Similarly, elements of $T(L(V, g))_{l}^{n}$ are finite sums of elements of the form 
$\varphi \psi$ where $\varphi\in T(L(V, g))_{l+a+i}$ 
and $\psi\in T(L(V, g))_{-\alpha-i}$ for $a\in \tilde{P}(V)$, $i\in n+\N$. Let $\eta\in T(L(V, g))_{k}$.
Then $\eta\varphi\in T(L(V, g))_{k+l+a+i}$ and 
\begin{eqnarray*}
\eta\varphi\psi &\in &T(L(V, g))_{k+l+a+i}T(L(V, g))_{-a-i}\nn
&\subset&T(L(V, g))_{k+l}^{n}\nn
&=&T(L(V, g))_{m}^{n}.
\end{eqnarray*}
Thus 
\begin{equation}\label{cont-3}
T(L(V, g))_{k}T(L(V, g))_{l}^{n_{2}}\subset T(L(V, g))_{m}^{n}.
\end{equation}

Combining (\ref{cont-2}) and (\ref{cont-3}), we obtain (\ref{cont-1}). 
\epfv

\begin{prop}\label{completion1}
For $k\in \tilde{P}(V)+\Z$,   elements of $\tilde{T}(L(V, g))_k$ are of the form $u_{1}+u_{2}$ where $u_{1}\in T(L(V, g))_k$ 
and $u_{2}$ is in the topological completion of 
$$\coprod_{a\in \tilde{P}(V), i\in \Z_{+}}T(L(V,g))_{k+a+i}T(L(V,g))_{-a-i}.$$
\end{prop}
\pf
The topological completion $\tilde{T}(L(V, g))_k$ of $T(L(V, g))_k$ is the space of all Cauchy sequences
or equivalently all Cauchy series in $T(L(V, g))_k$. Assuming that $\sum_{m\in \Z_{+}}
w_{m}$ where $w_{m}\in T(L(V, g))_{k}$ is such a Cauchy series. Then for any $n\in \Z$, 
$$\sum_{m=M}^{M+N} w_{m}\in \coprod_{a\in \tilde{P}(V), i\in n+\N}T(L(V,g))_{k+a+i}T(L(V,g))_{-a-i}$$
for sufficiently large $M\in \Z_{+}$ and $N\in \Z_{+}$. In particular, if we take $n=1$, then there exists
$M\in \Z_{+}$ such that $w_{m}\in  \coprod_{a\in \tilde{P}(V), i\in \Z_{+}}T(L(V,g))_{k+a+i}T(L(V,g))_{-a-i}$
for $m>M$. Let $u_{1}=\sum_{ m=1}^{M}w_{m}$ and $u_{2}=\sum_{ m=M+1}^{\infty}w_{m}$.
Then $u_{1}$ as a finite sum of elements of $T(L(V, g))_k$ is still an element of $T(L(V, g))_k$.
Since $\sum_{m\in \Z_{+}}
w_{m}$ is a Cauchy series, $u_{2}=\sum_{ m=M+1}^{\infty}w_{m}$ is also a Cauchy series. 
Since $w_{m}\in \coprod_{a\in \tilde{P}(V), i\in \Z_{+}}T(L(V,g))_{k+a+i}T(L(V,g))_{-a-i}$ when $m>M$,
it is in fact a Cauchy series in $\coprod_{a\in \tilde{P}(V), i\in \Z_{+}}T(L(V,g))_{k+a+i}T(L(V,g))_{-a-i}$
and thus is in its topological completion.
\epfv

Note that the two terms in the left-hand side of (\ref{jacobi-comp})  and the 
right-hand side of (\ref{jacobi-comp}) both correspond to well defined elements of $\tilde{T}(L(V, g))$.
Taking the difference between the sum of the two elements of $\tilde{T}(L(V, g))$  corresponding to the 
two terms in the left-hand side of (\ref{jacobi-comp}) and the element of $\tilde{T}(L(V, g))$ corresponding to 
the 
right-hand side of (\ref{jacobi-comp}), we obtain an element of $\tilde{T}(L(V, g))$. These are the elements 
of $\tilde{T}(L(V, g))$ corresponding to the Jacobi identity  (\ref{jacobi-comp}).
Let $J$ be the two-sided ideal in $\tilde{T}(L(V, g))$ generated by these elements of $\tilde{T}(L(V, g))$
corresponding to the Jacobi identity  (\ref{jacobi-comp}), the 
elements $\one(n)-\delta_{n, -1}$ for $n\in \Z$ corresponding to the identity property for $Y^{g}_{0}$,
and the elements 
\begin{equation}\label{L(-1)-der-comp}
(L(-1)u)(n) + (n+\mathcal{N})u(n-1)
\end{equation}
for $n\in \Z$ and $u\in V$ corresponding to the $L(-1)$-derivative property (\ref{L(-1)derivative}).

\begin{defn}
{\rm The {\it $g$-twisted universal enveloping algebra $U_g(V)$ of $V$} is  the quotient algebra 
$\tilde{T}(L(V, g))/J$ of $\tilde{T}(L(V, g))$ by the two-sided ideal $J$.
Since $J$ is generated by homogeneous elements, the grading on $\tilde{T}(L(V, g))$ induces
a grading on $U_g(V)$ and we still call the degree of an element of $U_g(V)$ its {\it weight}. The homogeneous subspace of 
$U_g(V)$ of weight $n\in \C$ is denoted by $U_{g}(V)_{n}$.}
\end{defn}

In the special case that $g=1_{V}$, the $g$-twisted universal enveloping algebra $U_g(V)$ of $V$
is the universal enveloping algebra $U(V)$ of $V$ introduced by Frenkel and Zhu \cite{FZ}.

The action of $g$ on $V^{[\alpha]}$ for $\alpha\in \C/\Z$ induces an action of $g$ on $L(V, g)$.
This action of $g$ on $L(V, g)$ induces an action of $g$ on $T(L(V, g))$ such that this action in fact 
gives an automorphism of the associative algebra $T(L(V, g))$. By the definition of the 
topological completion of $T(L(V, g))$, this action of $g$ on 
$T(L(V, g))$ extends continuously to an automorphism of the topological completion $\tilde{T}(L(V, g))$. 
Since the elements of
 $\tilde{T}(L(V, g))$ corresponding to 
(\ref{jacobi-comp}), the elements $\one(n)-\delta_{n, -1}$ and
(\ref{L(-1)-der-comp}) for $n\in \Z$  and $u$ in generalized eigenspaces for $g$ are also in generalized 
eigenspaces for $g$ in $\tilde{T}(L(V, g))$, $U_g(V)$ as the quotient of $\tilde{T}(L(V, g))$ by the ideal generated by 
these elements also has an action of $g$.  Moreover, the action of $g$ on $U_g(V)$  is an automorphism of $U_g(V)$.

The following result follows immediately from Proposition \ref{completion1}:

\begin{prop}\label{completion2}
For $k\in \C$,  let $\pi_{U_{g}(V)_{k}}$ be the projections from $\tilde{T}(L(V, g))_{k}$ to $U_{g}(V)_{k}$.
Then elements of $U_{g}(V)_{k}$ are of the form $u_{1}+u_{2}$ where 
$u_{1}\in \pi_{U_{g}(V)}T(L(V, g))_{k}$, and $u_{2}$ is  in the topological completion of 
$$\coprod_{a\in \tilde{P}(V), i\in \Z_{+}}(\pi_{U_{g}(V)}T(L(V,g))_{k+a+i})(\pi_{U_{g}(V)}T(L(V,g))_{-a-i}).\epfe$$ 
\end{prop}

\begin{prop}\label{ugv-mod}
Let $W$ be a $\overline{\C}_{+}$-graded weak $V$-module. Then $W$ has a natural structure of a $U_g(V)$-module such that 
the $U_g(V)$-module structure is compatible with the action of $g$ on $U_g(V)$ and $W$ in the 
sense that  for $u\in U_g(V)$ and $w\in W$,
$g(uw)=(gu)(gw)$. 
Conversely, let $W=\coprod_{n\in \overline{\C}_{+}, \alpha\in \C/\Z}W_{[n]}^{[\alpha]}$ be a $U_g(V)$-module equipped with 
gradings by $\overline{\C}_{+}$ and $\C/\Z$ and an action of $g$ on $W$
satisfying the following conditions:
\begin{enumerate}

\item The $U_g(V)$-module structure is compatible with the action of $g$ on $U_g(V)$ and $W$ in the sense that 
for $u\in U_g(V)$ and $w\in W$,
$g(uw)=(gu)(gw)$. 

\item For $u\in V^{[a+\Z]}$ $(a\in P(V))$ and $m\in \Z$ and $w\in W$,
$u(a+m)w=0$ when $m$ is sufficiently large.

\item For $u\in V_{(m)}$
 and $w\in W_{p}=\coprod_{\alpha\in \C/\Z}W_{p}^{[\alpha]}$ where 
$m\in \Z$ and $p\in \overline{\C}_{+}$,
$u(n)w$ is $0$ when $m-n-1+p\not\in \overline{\C}_{+}$,
is  in $W_{m-n-1+p}$ when $m-n-1+p\in \overline{\C}_{+}$ and for $r\in \R$,
$\coprod_{r\in \R}W_{\sqrt{-1}r}$ is equal to the subspace of $W$ consisting 
of $w\in W$ such that when  $u(n)w=0$ for 
$u\in V_{(m)}$, $m\in \Z$, $n\in \C$, $m-n-1 \not\in \overline{\C}_{+}$. 
For $\alpha \in \mathbb{C}/\mathbb{Z}$, $w \in
W^{[\alpha]}=\coprod_{n\in \overline{\C}_{+}}W_{n}^{[\alpha]}$, there exists $\Lambda \in \mathbb{Z}_{+}$
such that $ (g -
e^{2\pi\sqrt{-1}\alpha})^{\Lambda}w = 0$, where $L_{g}(0)=\omega(1)\in U_g(V)$.

\end{enumerate}
Then $W$ is a $\overline{\C}_{+}$-graded weak $g$-twisted $V$-module. In particular, the category of 
$\overline{\C}_{+}$-graded weak $g$-twisted $V$-modules and the category of $U_g(V)$-modules satisfying three
conditions above are isomorphic (not just equivalent since the underlying vector spaces are the same). 
Under this isomorphism, lower bounded generalized $g$-twisted $V$-modules
correspond to $U_g(V)$-modules of the form $W=\coprod_{n\in \C, \alpha\in \C/\Z}W_{[n]}^{[\alpha]}$
such that $W_{[n]}=\coprod_{\alpha\in \C/Z}W_{[n]}^{[\alpha]}$ for $n\in \C$ are generalized eigenspaces of 
$L_{g}(0)$, satisfying the three conditions above and the lower bounded condition that
for each $n\in \C$ and each $\alpha\in \C/\Z$, $W^{[\alpha]}_{[n + l]} = 0$ for
sufficiently negative real number $l$. Under this isomorphism, grading-restricted generalized $g$-twisted $V$-modules
(or simply $g$-twisted $V$-modules as we have called them) correspond to 
$U_g(V)$-modules of the form above, satisfying the three conditions above, the lower bounded condition and the condition that  for each $n \in
\mathbb{C}$, $\alpha \in \mathbb{C}/\mathbb{Z}$, 
$\dim W^{[\alpha]}_{[n]} < \infty$.
\end{prop}
\pf
The first conclusion is clear from the construction of $U_g(V)$. 

For its converse, we define 
$Y_{0}^{g}(u, x)w=\sum_{n\in \alpha+\Z}u(n)wx^{-n-1}$ for $u\in V$ and $w\in W$ and 
$Y^{g}(u, x)w=Y_{0}^{g}(x^{\mathcal{N}}u, x)w$. For $u\in V^{[a]}$, $m\in a+\Z$ and $w\in W^{[b+\Z]}$, we have
$u(m)w\in W^{[a+b+\Z]}$ Thus $gY_{0}^{g}(u, x)w=Y_{0}^{g}(gu, x)gw$ and $gY^{g}(u, x)w=Y^{g}(gu, x)gw$.
Since $u(a+m)w=0$  when $m\in \Z$ is sufficiently large,
$Y_{0}^{g}(u, x)w$ is lower truncated. Since $W$ is a $U_g(V)$-module, the Jacobi identity (\ref{Jacobi}) and the $L(-1)$-derivative property
(\ref{L(-1)derivative}) hold. By Theorem \ref{jacobi-equiv}, the duality property for $Y^{g}$ holds. It is also clear that 
all the other axioms for a $g$-twisted generalized $V$-module hold. 

The other conclusions are all clear.
\epfv

It is clear that the subspace $U_g(V)_{\mathbb{I}}=\coprod_{n\in \mathbb{I}}U_g(V)_{n}$ spanned by elements of $U_g(V)$
whose weights are imaginary numbers  is a subalgebra of $U_g(V)$. Let
$N(U_g(V)_{\mathbb{I}})$ be the two-sided ideal in $U_{g}(V)_{\mathbb{I}}$ generated by 
$$\coprod_{n\in \mathbb{I}, a\in \tilde{P}(V),i \in \N, \Re(a+i)>0}U_g(V)_{n+a+i} U_g(V)_{-a-i}.$$

\begin{defn}
{\rm The {\it $g$-twisted zero-mode algebra associated to $V$} is the quotient algebra
\[
Z_{g}(V) = U_g(V)_{\mathbb{I}}/N(U_g(V)_{\mathbb{I}}).
\]}
\end{defn}

For weight homogeneous $v \in V^{[a+\Z]}$,  $a\in P(V)\cap \mathbb{I}$, 
we write $o_{g}(v)$ for $v(\wt v +a- 1) \in U_g(V)_{\mathbb{I}}$, denote by $[o_{g}(v)]$ the image of
$o_{g}(v)$ in the quotient $Z_{g}(V)$, and extend $[o_{g}(v)]$ to all $v \in V$ by linearity. 

From the definition of $Z_{g}(V)$ and Proposition \ref{completion2}, we have the following useful result:

\begin{prop}\label{elements-z-g}
The $g$-twisted  zero-mode algebra $Z_{g}(V)$ is spanned by elements of the form 
$$[o_{g}(u_{1})]\cdots [o_{g}(u_{k})]$$
 for $u_{1}, \dots, u_{k}\in V$ and $k\in \N$. 
\end{prop}
\pf
By Proposition \ref{completion2}, elements of  $U_g(V)_{n}$ for $n\in \mathbb{I}$ 
are of the form $u_{1}+u_{2}$ where $u_{1}\in \pi_{U_{g}(V)}T(L(V, g))_{n}$, and $u_{2}$ is  in the topological completion of 
$$\coprod_{a\in \tilde{P}(V), i\in \Z_{+}}(\pi_{U_{g}(V)}T(L(V,g))_{n+a+i})(\pi_{U_{g}(V)}T(L(V,g))_{-a-i}).$$
Since $(\pi_{U_{g}(V)}T(L(V,g))_{n+a+i})(\pi_{U_{g}(V)}T(L(V,g))_{-a-i})\subset U_g(V)_{n+a+i}U_g(V)_{-a-i}$,
we see that $u_{2}$ is in the topological completion of $\coprod_{a\in \tilde{P}(V), i\in \Z_{+}}U_g(V)_{n+a+i}U_g(V)_{-a-i}$. 
Since $\Re(a+i)>0$ for $a\in \tilde{P}(V)$ and $i\in \Z_{+}$, we have 
$\coprod_{a\in \tilde{P}(V), i\in \Z_{+}}U_g(V)_{n+a+i}U_g(V)_{-a-i}\subset N(U_g(V)_{\mathbb{I}})$. 
Thus $u_{2}\in N(U_g(V)_{\mathbb{I}})$ and elements of  $Z_{g}(V)$ are of the form
$u_{1}+N(U_g(V)_{\mathbb{I}})$ for $u_{1}\in \pi_{U_{g}(V)}T(L(V, g))_{n}$, $n\in \mathbb{I}$ . 
But $T(L(V, g))_{n}$ is spanned by elements of the form $u_{1}(m_{1})\cdots u_{k}(m_{k})$ for 
$k\in \N$, homogeneous $u_{1}\in V^{[a_{1}+\Z]}, \dots, u_{k}\in V^{[a_{k}+\Z]}$, $a_{1}, \dots, a_{k}\in P(V)$,
 $m_{1}\in a_{1}+\Z, \dots, m_{k}\in a_{k}+\Z$ such that $\sum_{i=1}^{k}(\wt u_{i}-m_{i}-1)=n$. 
For simplicity, we shall also use the same notation to denote the images of such elements in $U_g(V)_{n}$
under the map $\pi_{U_{g}(V)}$. Since the commutator formula holds for elements in $U_{g}(V)$, 
these elements in $U_g(V)_{n}$ are linear combinations of elements of the same form but satisfying 
the additional condition $\Re(\wt u_{i}-m_{1}-1)\ge \cdots \ge \Re(\wt u_{i}-m_{k}-1)$. 
If $\Re(\wt u_{i}-m_{i}-1)<0$ for some $i$,  then $\Re(\wt u_{k}-m_{k}-1)<0$ and $u_{1}(m_{1})\cdots u_{k}(m_{k})$
 is in $U_{g}(V)_{n+a_{k}+(-\swt u_{k}+m_{k}-a_{k}+1)}U_g(V)_{-a_{k}-(-\swt u_{k}+m_{k}-a_{k}+1)}$ satisfying 
$\Re(a_{k}+(-\wt u_{k}+m_{k}-a_{k}+1))=\Re(-\wt u_{k}+m_{k}+1)>0$. Hence $u_{1}(m_{1})\cdots u_{k}(m_{k})$
 is in $N(U_g(V)_{\mathbb{I}})$. Thus if $u_{1}(m_{1})\cdots u_{k}(m_{k})$ is not in $N(U_g(V)_{\mathbb{I}})$,
$\Re(\wt u_{i}-m_{i}-1)\ge 0$ for $i=1, \dots, k$. Since $\sum_{i=1}^{k}(\wt u_{i}-m_{i}-1)=n\in \mathbb{I}$,
we must have $\Re(\wt u_{i}-m_{i}-1)= 0$ or $\Re(m_{i})=\wt u_{i}-1$ for $i=1, \dots, k$. Since $m_{i}\in a_{i}+\Z$,
we see that $\Re(a_{i})=0$ and hence $a_{i}\in \mathbb{I}$. Moreover, $m_{i}=\wt u_{i}+a_{i}-1$. 
Thus 
$$u_{1}(m_{1})\cdots u_{k}(m_{k})+N(U_g(V)_{\mathbb{I}})=[o_{g}(u_{1})]\cdots [o_{g}(u_{k})].$$
\epfv

We shall discuss modules for $U_g(V)$  and for $Z_{g}(V)$ and their connections with suitable $g$-twisted 
$V$-modules  in Section 5. In Section 6, Proposition \ref{elements-z-g} is needed in the proof that $Z_{g}(V)$ and $A_{g}(V)($ 
(to be constructed in the next section)
are isomorphic. Then as a consequence of the isomorphism, we obtain in turn an improvement of Proposition \ref{elements-z-g}. 

\setcounter{equation}{0}
\section{A generalization of Zhu's algebra to the $g$-twisted case for  $g$ of infinite order}

In \cite{DLM}, Dong, Li and Mason introduced a generalization of Zhu's algebra
for suitable $g$-twisted $V$-modules when $g$ is of finite order. But their generalization does not generalize straightforwardly to 
the case of a general automorphism $g$, especially in the case that $g$ does not act
on $V$ semisimply. In this section, we give such a generalization. When $g$ does not act on $V$ semisimply, the definition of this associative algebra 
involves the operator $\mathcal{N}$ obtained from the  multiplicative Jordan-Chevalley decomposition
of $g$. Moreover, this associative algebra is not a quotient of the fixed point subalgebra; instead, it is a quotient of  the 
generalized eigenspace of eigenvalue $1$ of $g$. This means that when $g$ does not act
on $V$ semisimply (in this case the order of $g$ must be infinite), this generalization of Zhu's algebra 
might have some interesting new features. 

In \cite{H1}, the first author found a natural definition of Zhu's algebra in connection 
with the modular invariance. Besides its geometric meaning in connection with the 
modular invariance, the advantage of this definition is that one can separate 
vertex operators from other formal series or analytic functions in the proofs of the 
properties of the algebra. This advantage makes it particularly simple for our generalization 
such that all the proofs are identical to those in \cite{H1} except for some properties 
of the operator $\mathcal{N}$. Because of this reason, we give our generalization 
by generalizing the definition in \cite{H1}. Just as in \cite{H1}, the difference between 
the definition given in this section and a direct generalization of the definitions of Zhu and Dong-Li-Mason
is given by an isomorphism obtained from the conformal transformation $\frac{1}{2\pi i}\log (1+2\pi i z)$,
the geometry underlying vertex operators and the Virasoro 
operators.

For weight homogeneous $u\in V^{[a+\Z]}$, $a\in P(V)$, let $l(u)\in a+\Z$ such that $\Re(\wt u-n-1)<0$ for $n\in a+\Z$ and 
$n-a\ge l(u)-a$. 
Then $l(u)=\wt u+a$ when $a\in P(V)\cap \mathbb{I}$ and $l(u)=\wt u+a-1$ when $a\in P(V)\cap ((0, 1)+\mathbb{I})$.  We define a linear operator 
$\L$ on $V$ by $\L u=l(u)u$ for weight homogeneous $u\in V^{[a+\Z]}$ and  linearity.
We define a product $\bullet_{g}$ in $V$ as follows: For $u, v\in V$,
\begin{eqnarray*}
u\bullet_{g} v&=&\res_{y}y^{-1}Y\left((1+y)^{\L-L(0)+\mathcal{N}}u, \frac{1}{2\pi \sqrt{-1}}\log(1+y)\right)v\nn
&=&\res_{x}\frac{2\pi \sqrt{-1}e^{2\pi\sqrt{-1} x}}{e^{2\pi\sqrt{-1}x}-1}
Y(e^{2\pi \sqrt{-1} x(\L-L(0)+\mathcal{N})}u, x)v.
\end{eqnarray*}
Recall the operator $\A_{V}$ defined in Section 2. Let $\tilde{O}_{g}(V)$ be the subspace of
$V$ spanned elements of the form 
\begin{eqnarray*}
\lefteqn{\res_{y}y^{-n}Y\left(y^{L(0)-\L+\A_{V}}(1+y)^{\L-L(0)+\mathcal{N}}u, \frac{1}{2\pi\sqrt{-1}}\log(1+y)\right)v}\nn
&&=\res_{x}\frac{2\pi \sqrt{-1}e^{2\pi\sqrt{-1} x}}{(e^{2\pi \sqrt{-1}x}-1)^{n}}
Y\left(\left(\frac{1}{e^{2\pi\sqrt{-1}x}-1}\right)^{L(0)-\L+\A_{V}}e^{2\pi \sqrt{-1}x(\L-L(0)+\mathcal{N})}u, x\right)v
\end{eqnarray*}
for $n>1$, $u, v\in V$.  Note that by definition, when $u \in V^{[\alpha]}$, $a\not\in \mathbb{I}$, $u\bullet_{g} v\in \tilde{O}_{g}(V)$.
Also note that 
\begin{equation}\label{tilde-O-g}
\res_{y}y^{-n}Y\left((1+y)^{\L-L(0)+\mathcal{N}}u, \frac{1}{2\pi\sqrt{-1}}\log(1+y)\right)v\in \tilde{O}_{g}(V)
\end{equation}
for $n>1$ since for homogeneous $u$, $(L(0)-\L+\A_{V})u$ is either $0$ or $u$.

Let $\tilde{A}_{g}(V)=V/\tilde{O}_{g}(V)$.

\begin{thm}\label{tilde-A-g}
The product $\bullet_{g}$ in $V$ induces a product (denoted still by $\bullet_{g}$)
in $\tilde{A}_{g}(V)$ such that $\tilde{A}_{g}(V)$ together with this product and
the equivalence class
of the vacuum $\mathbf{1}\in V$ is an associative 
algebra with identity. Moreover,
for $u\in V^{[a+\Z]}$, $a\in P(V)\cap\mathbb{I}$, $v\in V$, $((L(-1)+2\pi \sqrt{-1}(a+\mathcal{N}))u)\bullet_{g} v\in \tilde{O}_{g}(V)$ and for $u\in V^{[a+\Z]}$,
$v\in V^{[b+\Z]}$, $a, b\in P(V)\cap \mathbb{I}$, 
\begin{eqnarray}\label{bullet-bracket}
u\bullet_{g} v- v\bullet_{g} u&\equiv& 2\pi \sqrt{-1}\res_{x}
Y(e^{2\pi \sqrt{-1}x(a+\mathcal{N})}u, x)v \mod \tilde{O}_{g}(V)\nn
&=&2\pi \sqrt{-1} \sum_{i\in \N}\frac{2\pi \sqrt{-1}}{i!}((\mathcal{N}+a)^{i}u)_{i}v.
\end{eqnarray}
 In particular,  $\omega\bullet_{g} v-v\bullet_{g}\omega\equiv (2\pi \sqrt{-1} )^{2}\mathcal{N}v\mod \tilde{O}_{g}(V)$
for $v\in V^{[b+\Z]}$, $b\in P(V)\cap \mathbb{I}$, and 
$\omega+\tilde{O}_{g}(V)$ brackets with $v+\tilde{O}_{g}(V)$ sufficiently many times 
 in $\tilde{A}_{g}(V)$
is $0$ for $v\in V$.
\end{thm}
\pf
The proof of the first part is the same as the proof of the corresponding part in Proposition 6.1 in \cite{H1} except that 
we need to use (\ref{e-N-Y}) to take care of the additional operator $e^{2\pi \sqrt{-1}x\mathcal{N}}$
appearing in the definition and to take care of those factors involving $a\in P(V)$
just as in \cite{DLM}. 
To convince the reader that the same proof together with (\ref{e-N-Y})
and with the adjustments involving $a\in P(V)$ indeed works, we give as examples
the proof that
\begin{eqnarray}\label{assoc0}
&{\displaystyle \left(\res_{x}\frac{2\pi \sqrt{-1}e^{2\pi \sqrt{-1} x}}{(e^{2\pi \sqrt{-1}x}-1)^{n}}
Y\left(\left(\frac{1}{e^{2\pi \sqrt{-1}x}-1}\right)^{L(0)-\L+\A_{V}}e^{2\pi \sqrt{-1} x(\L-L(0)+\mathcal{N})}u_{1}, x\right)u_{2}\right)\bullet_{g} u_{3}}&\nn
&{\displaystyle =\left(\res_{x}\frac{2\pi \sqrt{-1}e^{2\pi \sqrt{-1} a x}}{(e^{2\pi \sqrt{-1}x}-1)^{n}}
Y(e^{2\pi \sqrt{-1} x\mathcal{N}}u_{1}, x)u_{2}\right)\bullet_{g} u_{3}}&
\end{eqnarray}
is in $\tilde{O}_{g}(V)$ when  $u_{1}\in V^{[a+\Z]}$ ($a\ne 0$), $u_{2}, u_{3}\in V$ and $n>1$.

In the case that $u_{2}\in V^{[b+\Z]}$ and $b\ne 1-a$, 
$$\res_{x}\frac{2\pi \sqrt{-1}e^{2\pi \sqrt{-1} ax}}{(e^{2\pi \sqrt{-1}x}-1)^{n}}
Y(e^{2\pi \sqrt{-1} x\mathcal{N}}u_{1}, x)u_{2}\in V^{[a+b+\Z]}\ne V^{[0]}.$$
Thus by definition, (\ref{assoc0}) is in $\tilde{O}_{g}(V)$. The only case we need to 
consider is  $u_{2}\in V^{[1-a+\Z]}$.

Using the definition of $\bullet_{g}$, (\ref{e-N-Y}), the Jacobi identity for the vertex operator for $V$ and 
the property of the formal delta function, we modify the first step in (6.8) in \cite{H1} to obtain that 
(\ref{assoc0}) is equal to 
\begin{eqnarray*}
\lefteqn{\res_{x_{0}}\res_{x_{2}}
\frac{2\pi \sqrt{-1}e^{2\pi \sqrt{-1} a x_{0}}}{(e^{2\pi \sqrt{-1}x_{0}}-1)^{n}}
\frac{2\pi \sqrt{-1}e^{2\pi \sqrt{-1} x_{2}}}{e^{2\pi \sqrt{-1}x_{2}}-1}
Y(e^{2\pi \sqrt{-1}x_{2}\mathcal{N}}Y(e^{2\pi \sqrt{-1}x_{0}\mathcal{N}}u_{1}, x_{0})u_{2}, x_{2})u_{3}}\nn
&&=\res_{x_{0}}\res_{x_{2}}
\frac{2\pi \sqrt{-1}e^{2\pi \sqrt{-1} a x_{0}}}{(e^{2\pi \sqrt{-1}x_{0}}-1)^{n}}
\frac{2\pi \sqrt{-1}e^{2\pi \sqrt{-1} x_{2}}}{e^{2\pi \sqrt{-1}x_{2}}-1}\cdot\nn
&&\quad\quad\quad\quad\quad\quad\cdot 
 Y(Y(e^{2\pi \sqrt{-1} (x_{2}+x_{0})\mathcal{N}}u_{1}, x_{0})e^{2\pi \sqrt{-1}x_{2}\mathcal{N}}u_{2}, x_{2})u_{3}\nn
&&=\res_{x_{1}}\res_{x_{2}}\res_{x_{0}}
\frac{2\pi \sqrt{-1}e^{2\pi \sqrt{-1} a x_{0}}}{(e^{2\pi \sqrt{-1}x_{0}}-1)^{n}}
\frac{2\pi \sqrt{-1}e^{2\pi \sqrt{-1} x_{2}}}{e^{2\pi \sqrt{-1}x_{2}}-1}\cdot\nn
&&\quad\quad\quad\quad\quad\quad\cdot 
x_{0}^{-1}\delta\left(\frac{x_{1}-x_{2}}{x_{0}}\right)
Y(e^{2\pi \sqrt{-1} (x_{2}+x_{0})\mathcal{N}}u_{1}, x_{1})Y(e^{2\pi \sqrt{-1}x_{2}\mathcal{N}}u_{2}, x_{2})u_{3}\nn
&&\quad -\res_{x_{1}}\res_{x_{2}}\res_{x_{0}}
\frac{2\pi \sqrt{-1}e^{2\pi \sqrt{-1} a x_{0}}}{(e^{2\pi \sqrt{-1}x_{0}}-1)^{n}}
\frac{2\pi \sqrt{-1}e^{2\pi \sqrt{-1} x_{2}}}{e^{2\pi \sqrt{-1}x_{2}}-1}\cdot\nn
&&\quad\quad\quad\quad\quad\quad\cdot 
x_{0}^{-1}\delta\left(\frac{x_{2}-x_{1}}{-x_{0}}\right)
Y(e^{2\pi \sqrt{-1}x_{2}\mathcal{N}}u_{2}, x_{2})Y(e^{2\pi \sqrt{-1} (x_{2}+x_{0})\mathcal{N}}u_{1}, x_{1})u_{3}\nn
&&=\res_{x_{1}}\res_{x_{2}}\res_{x_{0}}
\frac{2\pi \sqrt{-1}e^{2\pi \sqrt{-1} a x_{0}}}{(e^{2\pi \sqrt{-1}x_{0}}-1)^{n}}
\frac{2\pi \sqrt{-1}e^{2\pi \sqrt{-1} x_{2}}}{e^{2\pi \sqrt{-1}x_{2}}-1}\cdot\nn
&&\quad\quad\quad\quad\quad\quad\cdot 
x_{0}^{-1}\delta\left(\frac{x_{1}-x_{2}}{x_{0}}\right)
Y(e^{2\pi \sqrt{-1} x_{1}\mathcal{N}}u_{1}, x_{1})Y(e^{2\pi \sqrt{-1}x_{2}\mathcal{N}}u_{2}, x_{2})u_{3}
\end{eqnarray*}
\begin{eqnarray}\label{assoc1}
&&\quad -\res_{x_{1}}\res_{x_{2}}\res_{x_{0}}
\frac{2\pi \sqrt{-1}e^{2\pi \sqrt{-1} a x_{0}}}{(e^{2\pi \sqrt{-1}x_{0}}-1)^{n}}
\frac{2\pi \sqrt{-1}e^{2\pi \sqrt{-1} x_{2}}}{e^{2\pi \sqrt{-1}x_{2}}-1}\cdot\nn
&&\quad\quad\quad\quad\quad\quad\cdot 
x_{0}^{-1}\delta\left(\frac{x_{2}-x_{1}}{-x_{0}}\right)
Y(e^{2\pi \sqrt{-1}x_{2}\mathcal{N}}u_{2}, x_{2})Y(e^{2\pi \sqrt{-1} x_{1}\mathcal{N}}u_{1}, x_{1})u_{3}.
\end{eqnarray}
Then using the same remaining steps in (6.8) in \cite{H1}, we see that the right-hand side of 
(\ref{assoc1}) is equal to 
\begin{eqnarray}\label{assoc2}
\lefteqn{\sum_{k\in \N}\binom{-n}{k}\res_{x_{1}}\res_{x_{2}}2\pi \sqrt{-1} e^{2\pi \sqrt{-1} a (x_{1}-x_{2})}
(e^{2\pi \sqrt{-1} x_{1}}-1)^{-n-k}e^{-2\pi \sqrt{-1}(-n-k)x_{2}}\cdot}\nn
&&\quad \quad \cdot (e^{-2\pi \sqrt{-1} x_{2}}-1)^{k}\frac{2\pi \sqrt{-1}e^{2\pi \sqrt{-1} x_{2}}}{e^{2\pi \sqrt{-1}x_{2}}-1}
Y(e^{2\pi \sqrt{-1}x_{2}\mathcal{N}}u_{1}, x_{1})Y(e^{2\pi \sqrt{-1}x_{2}\mathcal{N}}u_{2}, x_{2})u_{3}\nn
&&\quad -\sum_{k\in \N}\binom{-n}{k}\res_{x_{1}}\res_{x_{2}}2\pi \sqrt{-1} e^{-2\pi \sqrt{-1} a (x_{2}-x_{1})}
(e^{-2\pi \sqrt{-1} x_{2}}-1)^{-n-k}e^{2\pi \sqrt{-1}(-n-k)x_{1}}\cdot\nn
&&\quad \quad \cdot (e^{2\pi \sqrt{-1} x_{1}}-1)^{k}\frac{2\pi \sqrt{-1}e^{2\pi \sqrt{-1} x_{2}}}{e^{2\pi \sqrt{-1}x_{2}}-1}
Y(e^{2\pi \sqrt{-1}x_{2}\mathcal{N}}u_{2}, x_{2})Y(e^{2\pi \sqrt{-1}x_{2}\mathcal{N}}u_{1}, x_{1})u_{3}\nn
&&=\sum_{k\in \N}\binom{-n}{k}\res_{x_{2}} e^{-2\pi \sqrt{-1} a x_{2}}
e^{-2\pi \sqrt{-1}(-n-k)x_{2}}(e^{-2\pi \sqrt{-1} x_{2}}-1)^{k}\frac{2\pi \sqrt{-1}e^{2\pi \sqrt{-1} x_{2}}}{e^{2\pi \sqrt{-1}x_{2}}-1}\cdot\nn
&&\quad \quad \cdot 
\res_{x_{1}}\frac{2\pi \sqrt{-1}e^{2\pi \sqrt{-1} a x_{1}}}{(e^{2\pi \sqrt{-1} x_{1}}-1)^{n+k}}
Y(e^{2\pi \sqrt{-1}x_{2}\mathcal{N}}u_{1}, x_{1})Y(e^{2\pi \sqrt{-1}x_{2}\mathcal{N}}u_{2}, x_{2})u_{3}\nn
&&\quad -\sum_{k\in \N}\binom{-n}{k}\res_{x_{1}} e^{2\pi \sqrt{-1} a x_{1}}(-1)^{-n-k}
e^{2\pi \sqrt{-1}(-n-k)x_{1}}(e^{2\pi \sqrt{-1} x_{1}}-1)^{k}\cdot\nn
&&\quad \quad \cdot 
\res_{x_{2}}\frac{2\pi \sqrt{-1} e^{2\pi \sqrt{-1} (1-a) x_{2}}}{(e^{2\pi \sqrt{-1}x_{2}}-1)^{1+n+k}}e^{2\pi \sqrt{-1} (n+k) x_{2}}
Y(e^{2\pi \sqrt{-1}x_{2}\mathcal{N}}u_{2}, x_{2})Y(e^{2\pi \sqrt{-1}x_{2}\mathcal{N}}u_{1}, x_{1})u_{3}.\nn
\end{eqnarray}
The first term in the right-hand side of (\ref{assoc2}) is in $\tilde{O}_{g}(V)$. Since $u_{2}\in V^{[1-a+\Z]}\ne V^{[0]}$ and 
$e^{2\pi \sqrt{-1} (n+k) x_{2}}$ can be written as a polynomial in $e^{2\pi \sqrt{-1}x_{2}}-1$ with degree less than 
or equal to $n+k$, the second term in the right-hand side of (\ref{assoc2}) is also in $\tilde{O}_{g}(V)$.

Next we prove $((L(-1)+2\pi \sqrt{-1} (a+\mathcal{N}))u)\bullet_{g} v\in \tilde{O}_{g}(V)$  for $u\in V^{[a+\Z]}$,
$a\in P(V)\cap \mathbb{I}$ and $v\in V$.
Using the properties of the vertex operator map $Y$ and (\ref{bracket-N-L}), we have
\begin{eqnarray*}
\lefteqn{(L(-1)u)\bullet_{g} v}\nn
&&=\res_{x}\frac{2\pi \sqrt{-1}e^{2\pi \sqrt{-1}x}}{e^{2\pi \sqrt{-1}x}-1}
Y(e^{2\pi \sqrt{-1}x (a+\mathcal{N})}L(-1)u, x)v\nn
&&=\res_{x}\frac{2\pi \sqrt{-1}e^{2\pi \sqrt{-1}x}}{e^{2\pi \sqrt{-1}x}-1}
\frac{d}{dx}Y(e^{2\pi \sqrt{-1}x(a+\mathcal{N})}u, x)v\nn
&& \quad-\res_{x}\frac{2\pi \sqrt{-1}e^{2\pi \sqrt{-1}x}}{e^{2\pi \sqrt{-1}x}-1}
Y(e^{2\pi \sqrt{-1}x (a+\mathcal{N})}2\pi \sqrt{-1} (a+\mathcal{N})u, x)v\quad\quad\quad\quad\quad\quad\quad\quad
\end{eqnarray*}
\begin{eqnarray}\label{l-1-N}
&&=-\res_{x}\left(\frac{d}{dx}\frac{2\pi \sqrt{-1}e^{2\pi \sqrt{-1}x}}{e^{2\pi \sqrt{-1}x}-1}\right)
Y(e^{2\pi \sqrt{-1}x (a+\mathcal{N})}u, x)v-2\pi \sqrt{-1} ((a+\mathcal{N})u)\bullet v\nn
&&=\res_{x}\frac{(2\pi \sqrt{-1})^{2}e^{2\pi \sqrt{-1}x}}{(e^{2\pi \sqrt{-1}x}-1)^{2}}
Y(e^{2\pi \sqrt{-1}x \mathcal{N}}v, x)\one-2\pi \sqrt{-1} ((a+\mathcal{N})u)\bullet v. 
\end{eqnarray}
Since the first term in (\ref{l-1-N}) is in $\tilde{O}_{g}(V)$, we obtain 
$((L(-1)+2\pi \sqrt{-1} (a+\mathcal{N}))u)\bullet_{g} v\in \tilde{O}_{g}(V)$. 

Now we prove (\ref{bullet-bracket}). For $u\in V^{[a+\Z]}$,
$v\in V^{[b+\Z]}$, $a, b\in P(V)\cap \mathbb{I}$, 
\begin{eqnarray*}
\lefteqn{u\bullet_{g} v-2\pi \sqrt{-1}\res_{x}
Y(e^{2\pi \sqrt{-1}x(a+\mathcal{N})}u, x)v}\nn
&&=\res_{x}\frac{2\pi \sqrt{-1}e^{2\pi \sqrt{-1}x}}{e^{2\pi \sqrt{-1}x}-1}
Y(e^{2\pi \sqrt{-1}x(a+\mathcal{N})}u, x)v-2\pi \sqrt{-1}\res_{x}
Y(e^{2\pi \sqrt{-1}x(a+\mathcal{N})}u, x)v\nn
&&=-\res_{x}\frac{2\pi \sqrt{-1}e^{-2\pi \sqrt{-1}x}}{e^{-2\pi \sqrt{-1}x}-1}
Y(e^{2\pi \sqrt{-1}x(a+\mathcal{N})}u, x)v\nn
&&=-\res_{x}\frac{2\pi \sqrt{-1}e^{-2\pi \sqrt{-1}x}}{e^{-2\pi \sqrt{-1}x}-1}
e^{2\pi \sqrt{-1}x(a+b+\mathcal{N})}Y(u, x)e^{-2\pi \sqrt{-1}x(b+\mathcal{N})}v\nn
&&\equiv-\res_{x}\frac{2\pi \sqrt{-1}e^{-2\pi \sqrt{-1}x}}{e^{-2\pi \sqrt{-1}x}-1}
e^{-xL(-1)}Y(u, x)e^{-2\pi \sqrt{-1}x(b+\mathcal{N})}v \mod \tilde{O}_{g}(V)\nn
&&=-\res_{x}\frac{2\pi \sqrt{-1}e^{-2\pi \sqrt{-1}x}}{e^{-2\pi \sqrt{-1}x}-1}
Y(e^{-2\pi \sqrt{-1}x(b+\mathcal{N})}v, -x)u \nn
&&=\res_{y}\frac{2\pi \sqrt{-1}e^{2\pi \sqrt{-1}y}}{e^{2\pi \sqrt{-1}y}-1}
Y(e^{2\pi \sqrt{-1}y(b+\mathcal{N})}v, y)u \nn
&&=v\bullet_{g} u,
\end{eqnarray*}
where we have used (\ref{e-N-Y}), the fact that the coefficients of $Y(u, x)e^{-2\pi \sqrt{-1}x(b+\mathcal{N})}v$
is in $V^{[a+b+\Z]}$,  (\ref{l-1-N}) with $a$ replaced by $a+b$ and the skew-symmetry of $Y$. 

Take $u=\omega$ and note that $\mathcal{N}\omega=0$, $a=0$ and $\omega_{0}=L(-1)$. Then we obtain 
$\omega\bullet_{g} v-v\bullet_{g} \omega=L(-1)v\equiv -2\pi \sqrt{-1}\mathcal{N}v \mod \tilde{O}_{g}(V)$. 
Since $\mathcal{N}$ is nilpotent, we see that the last conclusion is true.
\epfv

We now derive a generalization of the algebra $A_{g}(V)$ introduced in \cite{DLM}. 
Recall the invertible operator 
$$\mathcal{U}(1) = (2\pi \sqrt{-1})^{L(0)}e^{-L_{+}(A)}$$
on $V$ introduced in \cite{H1}, where $L_{+}(A)=\sum_{n\in \Z_{+}}A_{j}L(j)$ and the coeffieients $A_{j}$
for $j\in \Z_{+}$ are 
uniquely determined by
$$
\frac{1}{2\pi \sqrt{-1}}(e^{2\pi \sqrt{-1}y}-1)=\left(\exp\left(-\sum_{j\in \mathbb{Z}_{+}}
A_{j}y^{j+1}\frac{\partial}{\partial y}\right)\right)y.
$$
Note that $\mathcal{N}$ commutes with all the Virasoro operators (see (\ref{bracket-N-L})). Thus (1.5) 
in \cite{H1} gives 
\begin{eqnarray}\label{bullet=*}
\lefteqn{\mathcal{U}(1)\res_{y}y^{-n}
Y\left((1+y)^{\mathcal{N}}u, \frac{1}{2\pi \sqrt{-1}}\log(1+y)\right)v}\nn
&&=\res_{y}y^{-n}
Y\left((1+y)^{L(0)+\mathcal{N}}\mathcal{U}(1)u, y\right)\mathcal{U}(1)v.
\end{eqnarray}
We define a product $*$ on $V$ by
$$
u *_{g} v=\res_{x}x^{-1}Y\left((1+x)^{\L+\mathcal{N}}u, x\right)v
$$
for $u, v\in V$.
Let  $O_{g}(V)$ be the subspace of
$V$ spanned by elements of the form 
$$\res_{x}x^{-n}Y\left(x^{L(0)-\L+\A_{V}}(1+x)^{\L+\mathcal{N}}u, x\right)v$$
for $n>1$, $u, v\in V$.
By (\ref{bullet=*}),  we obtain the following result:

\begin{cor}\label{A-g}
Let $A_{g}(V)=V/O_{g}(V)$. Then $A_{g}(V)$
is an associative algebra isomorphic to $\tilde{A}_{g}(V)$. The isomorphism
from $\tilde{A}_{g}(V)$ to $A_{g}(V)$ is given by the map induced from 
$\mathcal{U}(1)$. \epf
\end{cor}

Because of Corollary \ref{A-g}, the categories of $A_{g}(V)$-modules and 
$\tilde{A}_{g}(V)$-modules are equivalent. 

Since the map $\mathcal{U}(1)$ corresponds to the coordinate change from the 
standard coordinate on an annulus centered at $0$ in the sphere $\C\cup \{0\}$
to the standard coordinate 
on the corresponding parallelogram, the algebra $\tilde{A}_{g}(V)$ defined using 
this coordinate change is most suitable for the study of the modular invariance
(see \cite{H1}).

\section{Functors between module categories}

In this section, we construct functors between suitable categories of 
$Z_{g}(V)$-modules, $\tilde{A}_{g}(V)$-modules, $A_{g}(V)$-modules and 
$g$-twisted
$V$-modules. 

We are interested only in $Z_{g}(V)$-modules that have compatible actions of $g$ and can be written as 
direct sums of generalized eigenspaces of the actions of $g$. Recall the notation $[o_{g}(v)]$ for $v\in V$ in Section 3. 
We are also interested in such $Z_{g}(V)$-modules with another grading 
given by the generalized eigenspaces of  $[o_{g}(\omega)]\in Z_{g}(V)$.

\begin{defn}
{\rm Let $M$ be a $Z_{g}(V)$-module.
If $M=\coprod_{\alpha\in \C/\Z}M^{[\alpha]}$ where $M^{[\alpha]}$ for 
$\alpha\in \C/\Z$ are  the generalized eigenspaces of the action of $g$, we call $M$ a {\it $g$-graded $Z_{g}(V)$-module}.  
Let $M=\coprod_{n\in \C, \alpha\in \C/\Z}M_{[n]}^{[\alpha]}$ be a $Z_{g}(V)$-module with double gradings by $\C$ and by $\C/\Z$ 
(or equivalently, by a subset of $[0, 1)$). 
If $M_{[n]}=\coprod_{\alpha\in \C/\Z} M_{[n]}^{[\alpha]}$ for $n\in \C$
 are generalized eigenspaces of the action of  $[o_{g}(\omega)]$ on $M$ with eigenvalues $n$ (called {\it weights}) and
the $\C/\Z$-grading is given by the generalized eigenspaces of an action of $g$ compatible with the $Z_{g}(V)$-module structure,
then we call $M$ a {\it doubly-graded $Z_{g}(V)$-module} or simply {\it graded $Z_{g}(V)$-module}.   
Let $M=\coprod_{n\in \C, \alpha\in \C/\Z}M_{[n]}^{[\alpha]}$ be a graded $Z_{g}(V)$-module, 
if for each $n\in \C$ and each $\alpha\in \C/\Z$, $M_{[n+l]}^{[\alpha]}=0$ for sufficiently negative real number $l$, we call $M$ a 
{\it lower bounded $Z_{g}(V)$-module}. A  lower bounded $Z_{g}(V)$-module $M=\coprod_{n\in \C, \alpha\in \C/\Z}M_{[n]}^{[\alpha]}$  
is called {\it grading restricted} if $\dim M_{[n]}^{[\alpha]}<\infty$ for $n\in \C, \alpha\in \C/\Z$. }
\end{defn}

Note that in particular, $Z_{g}(V)$ 
itself is a $g$-graded $Z_{g}(V)$-module. 

Let $W=\coprod_{n\in \overline{\C}_{+},\alpha\in \C/\Z}W_{n}^{[\alpha]}$ be a $\overline{\C}_{+}$-graded weak $g$-twisted $V$-module.
Recall the set $P(W)$ of all $g$-weights of $W$. 
In particular, $W=\coprod_{n\in \overline{\C}_{+}, a\in P(W)}W_{n}^{[a]}$. It is clear that for $a\in P(V)$ and $b\in P(W)$, 
either $a+b\in P(W)$ or $a+b-1\in P(W)$. 
Let 
$$\Omega^{g}(W) = \{w \in W\ |\ u(k)w = 0\ {\rm for\ homogeneous}\ u\in V, 
\Re(\wt u - k - 1)< 0\}.$$
Then by the $\overline{\C}_{+}$-grading condition in Definition \ref{d} and (\ref{yg-y0}),
$\Omega^{g}(W) =\coprod_{n\in \mathbb{I}}W_{n}=\coprod_{n\in \mathbb{I}, a\in P(W)}W_{n}^{[a+\Z]}$.
In the case that $W$ is a lower bounded $g$-twisted $V$-module, 
the gradings on $W$ also induce gradings on $\Omega^{g}(W)$. In particular, we have 
$$\Omega^{g}(W)=\coprod_{n\in \C}(\Omega^{g}(W))_{[n]}=\coprod_{a\in P(W)}(\Omega^{g}(W))^{[a+\Z]}
=\coprod_{n\in \C, a\in P(W)}(\Omega^{g}(W))^{[a+\Z]}_{[n]}.$$

By Theorem \ref{ugv-mod}, $W$ is a $U_{g}(V)$-module and hence also a $U_{g}(V)_{\mathbb{I}}$-module. 
In particular, $\Omega^{g}(W)$ generates a $U_{g}(V)_{\mathbb{I}}$-submodule of $W$.
By the commutator formula (\ref{comm2}), $\Omega^{g}(W) $ is invariant under
the action of $v(\wt v+a-1)$ for weight homogeneous $v\in V^{[a+\Z]}$, $a\in P(V)\cap \mathbb{I}$. 
In particular, the restriction of $v(\wt v+a-1)$ to $\Omega^{g}(W)$ gives us
a linear operator $o_{W}(v)$ on $\Omega^{g}(W)$. Also by the commutator formula (\ref{comm2}),
$u_{1}(m_{1})\cdots u_{k}(m_{k})w=0$ for homogeneous $u_{i}\in V^{[a_{i}]}$, $w\in \Omega^{g}(W)$ and $m_{i}\in a_{i}+\Z$
satisfying $\Re(\wt u_{1}-m_{1}-1+\cdots +\wt u_{k}-m_{k}-1+\wt w)<0$. 
Thus $N(U_{g}(V)_{\mathbb{I}})$ acts on $\Omega^{g}(W)$ as $0$. 
In particular, $\Omega^{g}(W)$ is a $U_{g}(V)_{\mathbb{I}}$-module and moreover  
$\Omega^{g}(W)$ is a $Z_{g}(V)$-module with the action of $Z_{g}(V)$  given by $[o_{g}(v)]w=o_{W}(v)w$ 
for $v\in V$ and $w\in \Omega^{g}(W)$. Note that the weight of the operator $v(\wt v+a-1)$
for weight homogeneous $v\in V^{[a+\Z]}$ where $a\in P(V)\cap \mathbb{I}$ is $-a\in \mathbb{I}$.
Thus we have:

\begin{prop}\label{z-module}
The action of $Z_{g}(V)$ on $\Omega^{g}(W)$ given by $[o_{g}(v)]w=o_{W}(v)w$  gives $\Omega^{g}(W)$ a
$g$-graded $Z_{g}(V)$-module structure. When $W$ is a lower bounded $g$-twisted $V$-module, $\Omega^{g}(W)$ is
a lower bounded $Z_{g}(V)$-module. In this case, 
$\coprod_{n\in r+\mathbb{I}}(\Omega^{g}(W))_{[n]}$ for $r\in \R$ are $Z_{g}(V)$-submodules
of $\Omega^{g}(W)$. \epf
\end{prop}

\begin{thm}\label{functor-H}
The functor $\Omega_{g}$ from the category of $\overline{\C}_{+}$-graded weak $g$-twisted $V$-modules
 to the category of $g$-graded $Z_{g}(V)$-modules given by $W\mapsto
\Omega_{g}(W)$ has a right inverse, that is, there exists a  functor $H_{g}$ from the
category of $g$-graded $Z_{g}(V)$-modules to the category of 
$\overline{\C}_{+}$-graded weak $g$-twisted $V$-modules such that $\Omega_{g} \circ H_{g}=1$, where 
$1$ is the identity functors on the category of $g$-graded $Z_{g}(V)$-modules.
 Moreover we can find such 
$H_{g}$ such that for any $\overline{\C}_{+}$-graded weak $g$-twisted $V$-module $W$, there exists 
a natural surjective homomorphism of $\overline{\C}_{+}$-graded weak $g$-twisted $V$-modules
from $H_{g}(\Omega_{g}(W))$ to 
the $\overline{\C}_{+}$-graded weak $g$-twisted $V$-submodule of $W$
generated by $\Omega_{g}(W)$.  The restriction of $\Omega_{g}$ to the category of lower bounded
$g$-twisted $V$-modules (or (grading-restricted generalized) $g$-twisted $V$-modules)   is a functor 
from this subcategory to the category of lower bounded $Z_{g}(V)$-modules (or grading-restricted $Z_{g}(V)$-modules).
The restrictions of $H_{g}$  to the category of  lower bounded $Z_{g}(V)$-modules (or grading-restricted $Z_{g}(V)$-modules)  is 
the right inverse of the restriction of $\Omega_{g}$ to the category of 
lower bounded
$g$-twisted $V$-modules (or (grading-restricted generalized) $g$-twisted $V$-modules)  such that for any 
object $W$ in the category, there exists a natural surjective homomorphism in the category
from $H_{g}(\Omega_{g}(W))$ to the
lower bounded
$g$-twisted $V$-submodule (or the (grading-restricted generalized) $g$-twisted $V$-submodule) of $W$
generated by $\Omega_{g}(W)$.
\end{thm}
\pf
Let $M=\coprod_{\alpha\in \C/\Z}M^{[\alpha]}$ be a $g$-graded $Z_{g}(V)$-module.
We now construct a $U_g(V)$-module $H_{g}(M)=\coprod_{n\in \overline{\C}_{+}, 
\alpha\in \C/\Z}H_{g}(M)_{n}^{[\alpha]}$  equipped with 
gradings by $\overline{\C}_{+}$ and $\C/\Z$ and an action of $g$ on $W$ satisfying Conditions 1 to 3
in Theorem \ref{ugv-mod}. Then
by Theorem \ref{ugv-mod}, we see that $H_{g}(M)$ is a $\overline{\C}_{+}$-graded weak $g$-twisted $V$-module. 

From the definition of $U_{g}(V)$, $U_{g}(V)_{-}=\coprod_{a\in \tilde{P}(V), i\in \N, \Re(a+i)>0}U_{g}(V)_{-a-i}$ is a subalgebra of $U_{g}(V)$. 
Since $M$ is a $Z_{g}(V)$-module, it is also a $U_{g}(V)_{\mathbb{I}}$-module. We define the action of 
$U_{g}(V)_{-}$ on $M$ to be $0$. Then $M$ becomes a  $U_{g}(V)_{-}\oplus U_{g}(V)_{\mathbb{I}}$-module.
Let $H_{g}(M)=U_{g}(V)\otimes_{U_{g}(V)_{-}\oplus U_{g}(V)_{\mathbb{I}}}M$. Then $H_{g}(M)$ is  a 
$U_g(V)$-module. By Proposition \ref{completion2}, elements of $U_{g}(V)_{k}$  are of the form $u_{1}+u_{2}$ where 
$u_{1}\in \pi_{U_{g}(V)}T(L(V, g))_{k}$, and $u_{2}$ is  in the topological completion of 
$$\coprod_{a\in \tilde{P}(V), i\in \Z_{+}}(\pi_{U_{g}(V)}T(L(V,g))_{k+a+i})(\pi_{U_{g}(V)}T(L(V,g))_{-a-i}).$$ 
Since $(\pi_{U_{g}(V)}T(L(V,g))_{-a-i})\subset U_{g}(V)_{-a-i}\subset U_{g}(V)_{-}$ for $a\in \tilde{P}(V)$ and $i\in \N$ 
satisfying $\Re(a+i)>0$, it
acts on $M$ as $0$ and thus 
$u_{2}w=0$ for $w\in M$.   Since 
$u_{1}\in \pi_{U_{g}(V)}T(L(V, g))_{k}$, it is 
a linear combination of elements of the form 
$u_{1}(m_{1})\cdots u_{k}(m_{k})$ for $u_{i}\in V^{[a_{i}+\Z]}$, $a_{i}\in P(V)$, $m_{i}\in a_{i}+\Z$.
Thus  elements of $H_{g}(M)$ are finite linear combinations of elements of the form 
$u_{1}(m_{1})\cdots u_{k}(m_{k})w$
for $u_{i}\in V^{[a_{i}+\Z]}$, $a_{i}\in P(V)$, $m_{i}\in a_{i}+\Z$ and $w\in M$. 
Using the commutator formula repeatedly and the actions of $U_{g}(V)_{-}$ and $U_{g}(V)_{\mathbb{I}}$
on $M$, we see that $H_{g}(M)$ is spanned by elements of the same form 
for homogeneous $u_{i}\in V^{[a_{i}+\Z]}$, $a_{i}\in P(V)$, $m_{i}\in a_{i}+\Z$ and $w\in M$,
satisfying $\Re(\wt u_{i}-m_{i}-1)>0$. 

For an element of the form $u_{1}(m_{1})\cdots u_{k}(m_{k})$ for weight  homogeneous $u_{i}\in V^{[a_{i}+\Z]}$, 
$a_{i}\in P(V)$, $m_{i}\in a_{i}+\Z$, we define its $\overline{\C}_{+}$-degree to be
$\wt u_{1}-m_{1}-1+\cdots +\wt u_{k}-m_{k}-1$. This gives a $\overline{\C}_{+}$-grading on $H_{g}(M)$. 
The actions of $g$ on $U_{g}(V)$ and $M$ also induce an action of $g$ on $H_{g}(M)$.
Moreover, the $\C/\Z$-gradings by generalized eigenspaces for the actions of $g$ on $U_{g}(V)$  and $M$ give 
a $\C/\Z$-grading by generalized eigenspaces for the action of $g$ on $H_{g}(M)$. The $\C/\Z$-grading
of $H_{g}(M)$ can also be given explicitly by defining the $\C/\Z$ degree of 
$u_{1}(m_{1})\cdots u_{k}(m_{k})w$ to be $\sum_{i=1}^{k}a_{i}+b+\Z$ when $w\in W^{[b+\Z]}$.

Condition 1 in Theorem \ref{ugv-mod} is satisfied by $H_{g}(M)$ because of the definition of the action of $g$ 
on  $H_{g}(M)$. Condition 2 in Theorem \ref{ugv-mod} is also satisfied by $H_{g}(M)$ because
$u(a+m)(u_{1}(m_{1})\cdots u_{k}(m_{k})w)=0$ when $m>\Re(\wt u-a-1+\sum_{i=1}^{k}(\wt u_{i}-m_{i}-1))$.
It is clear that $u_{n}$ changes the $\overline{\C}_{+}$-degree by $\wt u-n-1$ and 
the degree $0$ homogeneous subspace of  $H_{g}(M)$ is $M$ which is also the subspace of 
$H_{g}(M)$ annihilated by $U_{g}(V)_{-}$. Moreover, $H_{g}(M)$ is a direct sum of 
generalized eigenspaces for the action of $g$ on $H_{g}(M)$. Thus Condition 3 is satisfied.

Since our construction is natural, we obtain a functor $H_{g}$ from the category of $Z_{g}(V)$-modules
to the category of  $\overline{\C}_{+}$-graded weak $g$-twisted $V$-modules. 

By definition, $M\subset \Omega_{g}(H_{g}(M))$. Since $H_{g}(M)$ is spanned by elements of the form
$u_{1}(m_{1})\cdots u_{k}(m_{k})w$
for homogeneous $u_{i}\in V^{[a_{i}+\Z]}$, $a_{i}\in P(V)$, $m_{i}\in a_{i}+\Z$ and $w\in M$ satisfying $\wt u_{i}-m_{i}-1>0$, we see that 
$\Omega_{g}(H_{g}(M))=M$. Thus we have $\Omega_{g} \circ H_{g}=1$. 

Now let $W$ be a  $\overline{\C}_{+}$-graded weak $g$-twisted $V$-module. Then the 
$\overline{\C}_{+}$-graded weak $g$-twisted $V$-submodule
of $W$ generated by $\Omega_{g}(W)$ is spanned by elements of the form $u_{1}(m_{1})\cdots u_{k}(m_{k})w$
for $u_{i}\in V^{[a_{i}+\Z]}$, $a_{i}\in P(V)$, $m_{i}\in a_{i}+\Z$ and $w\in \Omega_{g}(W)$. 
Define a linear map from $H_{g}(\Omega_{g}(W))$ to this submodule of $W$ by sending the elements of the same form
in  $H_{g}(\Omega_{g}(W))$ to these corresponding elements in this submodule. This linear map is well defined 
because the only relations among these elements in $H_{g}(\Omega_{g}(W))$ are given by the commutator formula,
the weak associativity, the $L(-1)$-derivative property, the properties for the vacuum and the Virasoro
relations and these relations are all satisfied by the corresponding elements in the submodule of $W$. By definition, 
this linear map is a surjective module map from $H_{g}(\Omega_{g}(W))$ to this  submodule of $W$. 

The other conclusions follows immediately.
\epfv

We now discuss $\tilde{A}_{g}(V)$- and $A_{g}(V)$-modules and   $\overline{\C}_{+}$-graded weak $g$-twisted  $V$-modules.
As in the case of $Z_{g}(V)$-modules, we also discuss only $\tilde{A}_{g}(V)$- and $A_{g}(V)$-modules with actions of 
$g$ and are direct sums of generalized eigenspaces of the actions of $g$. As in the case of $Z_{g}(V)$-modules,
we are also interested in $g$-graded $\tilde{A}_{g}(V)$-modules 
with another grading 
given by the generalized eigenspaces of  $\omega+\tilde{O}_{g}(V)\in \tilde{A}_{g}(V)$.

\begin{defn}
{\rm Let $M$ be an  $\tilde{A}_{g}(V)$- or $A_{g}(V)$-module. If  $M=\coprod_{\alpha\in \C/\Z}M^{[\alpha]}$ where $M^{[\alpha]}$ for 
$\alpha\in \C/\Z$ are  the generalized eigenspaces of the action of $g$, then we call $M$
a {\it $g$-graded $\tilde{A}_{g}(V)$- and $A_{g}(V)$-module}.  
Let $M=\coprod_{n\in \C, \alpha\in \C/\Z}M_{[n]}^{[\alpha]}$ be an $\tilde{A}_{g}(V)$-module with double gradings by $\C$ and by $\C/\Z$ 
(or equivalently, by a subset of $[0, 1)+\mathbb{I}$). If $M_{[n]}=\coprod_{\alpha\in \C/\Z} M_{[n]}^{[\alpha]}$ for $n\in \C$
 are generalized eigenspaces of the action of  $\omega+\tilde{O}_{g}(V)$ on $M$ with eigenvalues $n$ (called {\it weights}) and
the $\C/\Z$-grading is given by the generalized eigenspaces of an action of $g$ compatible with the $\tilde{A}_{g}(V)$-module structure,
then we call $M$ a {\it doubly-graded $\tilde{A}_{g}(V)$-module} or simply a {\it graded $\tilde{A}_{g}(V)$-module}. 
Let $M=\coprod_{n\in \C, \alpha\in \C/\Z}M_{[n]}^{[\alpha]}$ be a graded $\tilde{A}_{g}(V)$-module, 
if for each $n\in \C$ and each $\alpha\in \C/\Z$, $M_{[n+l]}^{[\alpha]}=0$ for sufficiently negative real number $l$, we call $M$ a 
{\it lower bounded $\tilde{A}_{g}(V)$-module}. A  lower bounded $\tilde{A}_{g}(V)$-module $M=\coprod_{n\in \C, \alpha\in \C/\Z}M_{[n]}^{[\alpha]}$  
is called {\it grading restricted} if $\dim M_{[n]}^{[\alpha]}<\infty$ for $n\in \C, \alpha\in \C/\Z$. 
The notions of  {\it graded  $A_{g}(V)$-module}, {\it lower bounded $A_{g}(V)$-module} and {\it grading-restricted $A_{g}(V)$-module}
are defined in the same way. }
\end{defn}

Note that in particular, $\tilde{A}_{g}(V)$ ($A_{g}(V)$)
itself is a $g$-graded $\tilde{A}_{g}(V)$-module ($A_{g}(V)$-module).

Let $W=\coprod_{n\in \overline{\C}_{+},\alpha\in \C/\Z}W_{[n]}^{[\alpha]}$ be a $\overline{\C}_{+}$-graded weak $g$-twisted $V$-module.
Define two linear maps $o_{W}$ and $\rho_{W}=o_{W}\circ \mathcal{U}(1)$ from $V$ to $\mathrm{End}\;\Omega^{g}(W)$
by $v\mapsto o_{W}(v)$ and $v\mapsto o_{W}(\mathcal{U}(1)v)$, respectively, for $v\in V$, where
$$o_{W}(v)w=\res_{x}x^{a-1} Y_{0}(x^{L(0)}v, x)w=v(\wt v+a-1)w$$
for $a\in P(V)\cap \mathbb{I}$, weight homogeneous $v\in V^{[a+\Z]}$, $w\in \Omega^{g}(W)$ and
$$o_{W}(v)w=0$$
when $a\not\in \mathbb{I}$, $v\in V^{[a+\Z]}$, $w\in \Omega^{g}(W)$.
By definition, $o_{W}$ and $\rho_{W}$ are determined by 
their restrictions to $\coprod_{a\in P(V)\cap \mathbb{I}}V^{[a+\Z]}$.

We have the following:

\begin{thm}\label{module}
The spaces $O_{g}(V)$ and $\tilde{O}_{g}(V)$ are in the kernels of $o_{W}$ and $\rho_{W}$, respectively,
and the maps given by $u+O_{g}(V)\mapsto
o_{W}(u)$ and $u+\tilde{O}_{g}(V)\mapsto
\rho_{W}(u)$ for $u\in V$ give $\Omega_{g}(W)$ $g$-graded $A_{g}(V)$- and $\tilde{A}_{g}(V)$-module
structures, respectively. When $W$ is a lower bounded $g$-twisted $V$-module, $\Omega^{g}(W)$ is
a lower bounded $A_{g}(V)$- or $\tilde{A}_{g}(V)$-module. In particular, $\coprod_{n\in r+\mathbb{I}}(\Omega^{g}(W))_{[n]}$ for $r\in \R$ are
$A_{g}(V)$- or $\tilde{A}_{g}(V)$-submodules of $\Omega^{g}(W)$ . 
\end{thm}
\pf
As in the proofs of Theorem \ref{tilde-A-g} and Corollary \ref{A-g}, the proofs of the first statement for $\tilde{O}_{g}(V)$
and $\tilde{A}_{g}(V)$ are the same as the corresponding proofs of the first statement in Proposition 6.4 in \cite{H1} except for the use 
of  (\ref{e-N-Y}) and for the adjustments involving $a\in P(V)$. 
The proofs for $O_{g}(V)$ and $A_{g}(V)$ can be obtained either from those for 
$\tilde{O}_{g}(V)$
and $\tilde{A}_{g}(V)$ using the map $\mathcal{U}(1)$ or by directly generalizing the proof in \cite{Z}
and \cite{DLM}. 

Here we give only the proof of the statement that 
\begin{equation}\label{module0}
\rho_{W}\left(\res_{x_{0}}
\frac{2\pi \sqrt{-1}e^{2\pi \sqrt{-1} a x_{0}}}{(e^{2\pi \sqrt{-1}x_{0}}-1)^{n}}Y(e^{2\pi \sqrt{-1}x_{0}\mathcal{N}}u, x_{0})v\right)=0
\end{equation}
for $u\in V^{[a+\Z]}$, $a\in P(V)$, $0<\Re(a)<1$,  and $v\in V^{[b+\Z]}$.

If $a+b-1\not\in \mathbb{I}$, $Y(e^{2\pi \sqrt{-1}x_{0}\mathcal{N}}u, x_{0})v
\in V^{[a+b+\Z]}[[x, x^{-1}]]$ where $a+b+\Z\not\not \in \mathbb{I}+\Z$. Then  
by definition, 
$$\rho_{W}(Y(e^{2\pi \sqrt{-1}x_{0}\mathcal{N}}u, x_{0})v)=0$$
 and thus 
(\ref{module0}) holds. So we need 
only consider the case $a+b-1\in \mathbb{I}$. In this case,
for $n\in \mathbb{Z}_{+}$ and $w\in \Omega_{g}(W)$,  
using the definitions, (1.5) in \cite{H1}, the 
$L(0)$-conjugation formula, (\ref{bracket-N-L}), the basic property of the formal delta-function 
and the Jacobi identity (\ref{Jacobi}) for $Y_{0}$, we have 
\begin{eqnarray*}
\lefteqn{\rho_{W}\left(\res_{x_{0}}
\frac{2\pi \sqrt{-1}e^{2\pi \sqrt{-1}a x_{0}}}{(e^{2\pi \sqrt{-1}x_{0}}-1)^{n}}Y(e^{2\pi \sqrt{-1}x_{0}\mathcal{N}}u, x_{0})v\right)w}\nn
&&=\res_{x_{2}}\frac{1}{x_{2}^{(1-a-b)+1}}\res_{x_{0}}
\frac{2\pi \sqrt{-1}e^{2\pi \sqrt{-1}a x_{0}}}{(e^{2\pi \sqrt{-1}x_{0}}-1)^{n}}
Y_{0}(\mathcal{U}(x_{2})Y(e^{2\pi \sqrt{-1}x_{0}\mathcal{N}}u, x_{0})v, x_{2})w\nn
\end{eqnarray*}
\begin{eqnarray}\label{module1}
&&=\res_{x_{2}}\res_{x_{0}}
\frac{2\pi \sqrt{-1}e^{2\pi \sqrt{-1}a x_{0}}}{x_{2}^{2-a-b}(e^{2\pi \sqrt{-1}x_{0}}-1)^{n}}\cdot \nn
&&\quad\quad\quad\cdot 
Y_{0}(x_{2}^{L(0)}Y(\mathcal{U}(e^{2\pi \sqrt{-1}x_{0}})e^{2\pi \sqrt{-1}x_{0}\mathcal{N}}u, e^{2\pi \sqrt{-1}x_{0}}-1)
\mathcal{U}(1)v, x_{2})w\nn
&&=\res_{x_{2}}\res_{x_{0}}
\frac{2\pi \sqrt{-1}e^{2\pi \sqrt{-1} a x_{0}}}{x_{2}^{2-a-b}(e^{2\pi \sqrt{-1}x_{0}}-1)^{n}}\cdot \nn
&&\quad\quad\quad\cdot Y_{0}(Y(x_{2}^{L(0)}\mathcal{U}(e^{2\pi \sqrt{-1}x_{0}})e^{2\pi \sqrt{-1}x_{0}\mathcal{N}}u, 
x_{2}(e^{2\pi \sqrt{-1}x_{0}}-1))
x_{2}^{L(0)}\mathcal{U}(1)v, x_{2})w\nn
&&=\res_{x_{2}}\res_{x_{0}}
\frac{2\pi \sqrt{-1}e^{2\pi \sqrt{-1}a x_{0}}}{x_{2}^{2-a-b}(e^{2\pi \sqrt{-1}x_{0}}-1)^{n}}\cdot \nn
&&\quad\quad\quad\cdot Y_{0}(Y(e^{2\pi \sqrt{-1}x_{0}\mathcal{N}}x_{2}^{L(0)}\mathcal{U}(e^{2\pi \sqrt{-1}x_{0}})u, 
x_{2}(e^{2\pi \sqrt{-1}x_{0}}-1))
x_{2}^{L(0)}\mathcal{U}(1)v, x_{2})w\nn
&&=\res_{x_{2}}\res_{y_{0}}\left(1+\frac{y_{0}}{x_{2}}\right)^{a-1}
\frac{1}{x_{2}^{3-a-b-n}y_{0}^{n}}\cdot\nn
&&\quad\quad\quad \cdot
Y_{0}\left(Y\left(\left(1+\frac{y_{0}}{x_{2}}\right)^{\mathcal{N}}\mathcal{U}(x_{2}+y_{0})u, y_{0}\right)
\mathcal{U}(x_{2})v, x_{2}\right)w\nn
&&=\res_{x_{2}}\res_{y_{0}}\res_{x_{1}}x_{1}^{a-1}x_{1}^{-1}\delta\left(\frac{x_{2}+y_{0}}{x_{1}}\right)
\left(\frac{x_{2}+y_{0}}{x_{1}}\right)^{a}
\frac{x_{2}^{b-2+n}}{y_{0}^{n}}\cdot\nn
&&\quad\quad\quad \cdot
Y_{0}\left(Y\left(\left(1+\frac{y_{0}}{x_{2}}\right)^{\mathcal{N}}\mathcal{U}(x_{2}+y_{0})u, y_{0}\right)
\mathcal{U}(x_{2})v, x_{2}\right)w\nn
&&=\res_{x_{2}}\res_{y_{0}}\res_{x_{1}}x_{1}^{a-1}x_{1}^{-1}\delta\left(\frac{x_{2}+y_{0}}{x_{1}}\right)
\left(\frac{x_{2}+y_{0}}{x_{1}}\right)^{a}
\frac{x_{2}^{b-2+n}}{y_{0}^{n}}\cdot\nn
&&\quad\quad\quad \cdot
Y_{0}\left(Y\left(\left(1+\frac{y_{0}}{x_{2}}\right)^{\mathcal{N}}\mathcal{U}(x_{1})u, y_{0}\right)
\mathcal{U}(x_{2})v, x_{2}\right)w\nn
&&=\res_{x_{2}}\res_{y_{0}}\res_{x_{1}}
\frac{x_{1}^{a-1}x_{2}^{b-2+n}}{y_{0}^{n}}y_{0}^{-1}\delta\left(\frac{x_{1}-x_{2}}{y_{0}}\right)
\cdot \nn
&&\quad\quad\quad\cdot Y_{0}(\mathcal{U}(x_{2}+y_{0})u, x_{1})
Y_{0}(\mathcal{U}(x_{2})v, x_{2})w\nn
&&\quad -\res_{x_{2}}\res_{y_{0}}\res_{x_{1}}
\frac{x_{1}^{a-1}x_{2}^{b-2+n}}{y_{0}^{n}}y_{0}^{-1}\delta\left(\frac{x_{2}-x_{1}}{-y_{0}}\right)
\cdot \nn
&&\quad\quad\quad\cdot Y_{0}(\mathcal{U}(x_{2})v, x_{2})
Y_{0}(\mathcal{U}(x_{2}+y_{0})u, x_{1})w\nn
&&=\res_{x_{2}}\res_{x_{1}}\frac{x_{1}^{a-1}x_{2}^{b-2+n}}{(x_{1}-x_{2})^{n}}
Y_{0}(\mathcal{U}(x_{1})u, x_{1})
Y_{0}(\mathcal{U}(x_{2})v, x_{2})w\nn
&&\quad -\res_{x_{2}}\res_{x_{1}}
\frac{x_{1}^{a-1}x_{2}^{b-2+n}}{(-x_{2}+x_{1})^{n}}
Y_{0}(\mathcal{U}(x_{2})v, x_{2})
Y_{0}(\mathcal{U}(x_{1})u, x_{1})w.
\end{eqnarray}
Since $w\in \Omega_{g}(W)$, the right-hand side of (\ref{module1}) is $0$.

Since the weight of the operator $\rho_{W}(u)$ is $\mathbb{I}$, $\coprod_{n\in r+\mathbb{I}}(\Omega^{g}(W))_{[n]}$ for $r\in \R$ are
invariant under  $\rho_{W}(u)$ and thus are $\tilde{A}(V)$-submodules of $\Omega^{g}(W)$. The compatibility 
with the action of $g$ clear.
\epfv

\begin{thm}\label{functor-f}
The functor $\Omega_{g}$ from the category of $\overline{\C}_{+}$-graded weak $g$-twisted $V$-modules
 to the category of  $g$-graded $A_{g}(V)$-modules or $\tilde{A}_{g}(V)$-modules given by $W\mapsto
\Omega_{g}(W)$ has a right inverse, that is, there exist functors $S_{g}$ and $\tilde{S}_{g}$ from the
categories of $g$-graded $A_{g}(V)$-modules  and $\tilde{A}_{g}(V)$-modules, respectively, to the category of 
$\overline{\C}_{+}$-graded weak $g$-twisted $V$-modules such that $\Omega_{g} \circ S_{g}=1$ and 
$\Omega_{g} \circ \tilde{S}_{g}=\tilde{1}$, where 
$1$ and $\tilde{1}$ are the identity functors on the categories of $A_{g}(V)$-modules and $\tilde{A}_{g}(V)$-modules, respectively.
 Moreover we can find such 
$S_{g}$ and $\tilde{S}_{g}$ such that for any 
$\overline{\C}_{+}$-graded weak $g$-twisted $V$-module
 $W$, there exists a natural surjective homomorphism of $\overline{\C}_{+}$-graded weak $g$-twisted $V$-modules
from $S_{g}(\Omega_{g}(W))$ or $\tilde{S}_{g}(\Omega_{g}(W))$ to 
the $\overline{\C}_{+}$-graded weak $g$-twisted $V$-submodule of $W$
generated by $\Omega_{g}(W)$.  The restriction of $\Omega_{g}$ to the category of lower bounded
$g$-twisted $V$-modules (or (grading-restricted) $g$-twisted $V$-mdoules) is a functor 
from this subcategory to the category of lower bounded (or grading-restricted) $A_{g}(V)$-modules or $\tilde{A}_{g}(V)$-modules.
The restrictions of $S_{g}$ and $\tilde{S}_{g}$ to the categories of lower bounded (or grading-restricted) $A_{g}(V)$-modules and $\tilde{A}_{g}(V)$-modules,
respectively, are right inverses of the restriction of $\Omega_{g}$ to the categories of 
lower bounded (or grading-restricted) $g$-twisted $V$-modules such that for any lower bounded (or grading-restricted) 
 $g$-twisted
$V$-module $W$, there exists a natural surjective homomorphism of such $g$-twisted $V$-modules
from $S_{g}(\Omega_{g}(W))$ or $\tilde{S}_{g}(\Omega_{g}(W))$ to 
the lower bounded (or grading-restricted) $g$-twisted $V$-submodule of $W$
generated by $\Omega_{g}(W)$.
\end{thm}
\pf
We prove only the results for $\tilde{A}_{g}(V)$. The proofs of the results for $A_{g}(V)$ can be derived 
from the results for $\tilde{A}_{g}(V)$.
The proofs of these results are the same as the proofs of the second and third statements of Proposition 6.4 in \cite{H1}
except that instead of expressions corresponding to the commutator formula, the $L(-1)$-derivative property and 
the associator formula
for $V$-modules, here we use the expressions corresponding to (\ref{commu}), (\ref{L(-1)derivative}) and (\ref{asso}). 

We need to construct a
$\overline{\C}_{+}$-graded weak $g$-twisted $V$-module from a $g$-graded $\tilde{A}_{g}(V)$-module.
Recall the space $L(V, g)$, the tensor algebra $T(L(V, g))$ and the notations we have used 
in Section 3. 

Let $M=\coprod_{\alpha\in \C/\Z}M^{[\alpha]}$ be a $g$-graded $\tilde{A}_{g}(V)$-module and let $\rho: \tilde{A}_{g}(V)
\to \edo M$ be the map giving the representation of $\tilde{A}_{g}(V)$
on $M$. Consider 
$T(L(V, g))\otimes M$. Again for simplicity we shall omit the 
tensor product sign. So $T(L(V, g))\otimes M$
is spanned by elements
of the form $u_{1}(m_{1})\cdots u_{k}(m_{k})w$ for
$u_{i}\in V^{[a_{i}+\Z]}$, $m_{i}\in a_{i}+\mathbb{Z}$, $a_{i}\in P(V)$,  $i=1, \dots, k$, and $w\in M^{[b+\Z]}_{[n]}$,
$b\in [0, 1)+\mathbb{I}$, $n\in \C$. 
For $u\in V^{[a+\Z]}$, $m\in a+\mathbb{Z}$ where $a\in P(V)$, $u(m)$ acts from the left on 
$T(L(V, g))\otimes M$.
In the case that $u_{i}$ are homogeneous with respect to the grading by weights, 
we define the $\overline{\C}_{+}$-degree  of $u_{1}(m_{1})\cdots u_{k}(m_{k})w$
to be $(\wt u_{1}-m_{1}-1)+\cdots +(\wt u_{k}-m_{k}-1)$ and its $\C/\Z$ degree by 
$a_{1}+\cdots +a_{k}+b+\Z$. 
For any $u\in V^{[a+\Z]}$, let
$$Y_{0}^{t}(u, x): T(L(V, g))\otimes M
\to T(L(V, g))[[x, x^{-1}]]$$
be defined by 
$$Y_{0}^{t}(u, x)=\sum_{m\in a+\mathbb{Z}}u(m)x^{-m-1}.$$
For a homogeneous element $u\in V^{[a+\Z]}$, $a\in P(V)\cap \mathbb{I}$, let $o_{t}(u)=u(\wt u+a-1)$. 
Using linearity, we extend $o_{t}(u)$ to inhomogeneous $u\in \coprod_{a\in P(V)\cap \mathbb{I}}V^{[a+\Z]}$.

Let $\mathcal{I}$ be the graded 
$T(L(V, g))$-submodule 
of $T(L(V, g))\otimes M$
generated by elements of the forms
$u(m)w$ ($u\in V^{[a+\Z]}$, $a\in P(V)$, $m\in a+\Z$, $\wt u-m-1< 0$, $w\in M$), 
$o_{t}(\mathcal{U}(1)u)w-\rho(u+\tilde{O}(V))w$ ($u\in \coprod_{a\in P(V)\cap \mathbb{I}}V^{[a+\Z]}$, $w\in M$) and
the coefficients in $x_{1}$ and $x_{2}$ of 
\begin{eqnarray}\label{commu-rel}
\lefteqn{Y_{0}^{t}(u, x_{1})Y_{0}^{t}(v, x_{2})w
-Y_{0}^{t}(v, x_{2})Y_{0}^{t}(u, x_{1})w}\nn
&&-\res_{x_{0}}x_{1}^{-1}\delta\left(\frac{x_{2}+x_{0}}{x_{1}}\right)\left(\frac{x_{2}+x_{0}}{x_{1}}\right)^{a}
Y_{0}^{t}\left(Y\left(\left(1+\frac{x_{0}}{x_{2}}\right)^{\mathcal{N}}u, x_{0}\right)v, x_{2}\right)w
\end{eqnarray}
($u\in V^{[a+\Z]}, v\in V$, $a\in P(V)$ and $w\in T(L(V, g))\otimes M$). 
Note that as in the proof of Proposition 6.4 in \cite{H1}, the coefficients of the formal expression above 
are indeed in $T(L(V, g))\otimes M$.
Let $S_{g}^{1}(M)=(T(L(V, g))\otimes M)/\mathcal{I}$. Then 
$S_{g}^{1}(M)$ is also  a graded 
$T(L(V, g))$-module. 
In fact, by the definition of $\mathcal{I}$, we see that 
$S_{g}^{1}(M)$ is  spanned by elements
of the form $u_{1}(m_{1})\cdots u_{k}(m_{k})w+\mathcal{I}$ for $k\in \N$, homogeneous
$u_{i}\in V^{[a_{i}+\Z]}$, $m_{i}\in a_{i}+\mathbb{Z}$, $a_{i}\in P(V)$,  $i=1, \dots, k$, and $w\in M^{[b+\Z]}$,
$b\in [0, 1)+\mathbb{I}$, satisfying 
$\Re(\wt u_{1}-m_{1}-1)\ge \cdots\ge \Re(\wt u_{k}-m_{k}-1)>0$. In particular,  $S_{g}^{1}(M)$ is doubly graded by $\overline{\C}_{+}$ and 
$\C/\Z$.
Thus for homogeneous $u\in V$ and $w\in S_{g}^{1}(M)$,
$u(m)w=0$ when the real part of $m$ is sufficiently large. Moreover, $S_{g}^{1}(M)=S_{g}^{1}(M)_{+}\oplus S_{g}^{1}(M)_{0}$, where 
$S_{g}^{1}(M)_{+}$ is the subspace of $S_{g}^{1}(M)$ spanned by elements of the form above 
for $k\in \Z_{+}$ and $S_{g}^{1}(M)_{0}$ is spanned by the elements of the form $w+\mathcal{I}$ for $w\in M$. 

We show that 
$S_{g}^{1}(M)_{0}$ is linearly isomorphic to $M$. We need only show that if $\sum_{i=1}^{k}\rho(u^{i}_{1}+\tilde{O}(V))
\cdots \rho(u^{i}_{k}+\tilde{O}(V))w\ne 0$ in $M$, then 
$\sum_{i=1}^{k}\rho(u^{i}_{1}+\tilde{O}(V))
\cdots \rho(u^{i}_{k}+\tilde{O}(V))w+\mathcal{I}\ne 0$ in $S_{g}^{1}(M)_{0}$, or equivalently, 
$\sum_{i=1}^{k}o_{t}(\mathcal{U}(1)u^{i}_{1})
\cdots o_{t}(\mathcal{U}(1)u^{i}_{k})w\not\in \mathcal{I}$. 
Assume that  this is not true, that is, $\sum_{i=1}^{k}o_{t}(\mathcal{U}(1)u^{i}_{1})
\cdots o_{t}(\mathcal{U}(1)u^{i}_{k})w\in \mathcal{I}$. Since $\sum_{i=1}^{k}\rho(u^{i}_{1}+\tilde{O}(V))
\cdots \rho(u^{i}_{k}+\tilde{O}(V))w\ne 0$ in $M$, $\sum_{i=1}^{k}o_{t}(\mathcal{U}(1)u^{i}_{1})
\cdots o_{t}(\mathcal{U}(1)u^{i}_{k})w$ cannot be in the  graded 
$T(L(V, g))$-submodule 
of $T(L(V, g))\otimes M$
generated by elements of the forms
$o_{t}(\mathcal{U}(1)u)w-\rho(u+\tilde{O}(V))w$ ($u\in \coprod_{a\in P(V)\cap \mathbb{I}}V^{[a+\Z]}$, $w\in M$).
Clearly, $\sum_{i=1}^{k}o_{t}(\mathcal{U}(1)u^{i}_{1})
\cdots o_{t}(\mathcal{U}(1)u^{i}_{k})w$ is not in the  graded 
$T(L(V, g))$-submodule 
of $T(L(V, g))\otimes M$
generated by elements of the forms $u(m)w$ ($u\in V^{[a+\Z]}$, $a\in P(V)$, $m\in a+\Z$, $\wt u-m-1< 0$, $w\in M$).
Thus $\sum_{i=1}^{k}o_{t}(\mathcal{U}(1)u^{i}_{1})
\cdots o_{t}(\mathcal{U}(1)u^{i}_{k})w$ is in fact in the 
the graded 
$T(L(V, g))$-submodule $\widetilde{\mathcal{I}}$
of $T(L(V, g))\otimes M$
generated by 
the coefficients in $x_{1}$ and $x_{2}$ of (\ref{commu-rel}) ($u\in V^{[a+\Z]}, v\in V$, $a\in P(V)$ and $w\in T(L(V, g))\otimes M$). 
Then by the definition of $\widetilde{\mathcal{I}}$, 
$\sum_{i=1}^{k}o_{t}(\mathcal{U}(1)u_{1})
\cdots o_{t}(\mathcal{U}(1)u_{k})$ must be in the subalgebra of $T(L(V, g))$ spanned by products of 
at least one element of the form 
$$o_{t}(\mathcal{U}(1)v_{1})o_{t}(\mathcal{U}(1)v_{2})-o_{t}(\mathcal{U}(1)v_{2})o_{t}(\mathcal{U}(1)v_{1})
-2\pi \sqrt{-1}\res_{x_{0}}o_{t}(\mathcal{U}(1)Y(e^{2\pi \sqrt{-1}x(a_{1}+\mathcal{N})}v_{1}, x_{0})v_{2})$$
for $v_{j}\in V^{[a_{j}+\Z]}$, $a_{j}\in P(V)\cap \mathbb{I}$ for $j=1, 2$, and elements of the form 
$o_{t}(\mathcal{U}(1)v)$ for $v\in V^{[a+\Z]}$, $a\in P(V)\cap \mathbb{I}$. So the element 
$\sum_{i=1}^{k}\rho(u^{i}_{1}+\tilde{O}(V))
\cdots \rho(u^{i}_{k}+\tilde{O}(V))w$ of $M$ corresponding to 
$\sum_{i=1}^{k}o_{t}(\mathcal{U}(1)u^{i}_{1})
\cdots o_{t}(\mathcal{U}(1)u^{i}_{k})w$ must be a sum of elements obtained by 
applying at least one operator of the form 
\begin{eqnarray*}
\lefteqn{\rho(v_{1}+\tilde{O}(V))\rho(v_{2}+\tilde{O}(V))-\rho(v_{2}+\tilde{O}(V))\rho(v_{1}+\tilde{O}(V))}\nn
&& \quad -2\pi \sqrt{-1}\res_{x_{0}}\rho(Y(e^{2\pi \sqrt{-1}x(a_{1}+\mathcal{N})}v_{1}, x_{0})v_{2}+\tilde{O}(V))\nn
&&=\rho(v_{1}\bullet_{g}v_{2}-v_{2}\bullet_{g} v_{1}-2\pi \sqrt{-1}\res_{x_{0}}Y(e^{2\pi \sqrt{-1}x(a_{1}+\mathcal{N})}v_{1}, x_{0})v_{2}+\tilde{O}(V))
\end{eqnarray*}
and operators of the form $\rho(v)$ for $v\in V^{[a+\Z]}$, $a\in P(V)\cap \mathbb{I}$ to $w$.
By (\ref{bullet-bracket}), 
$$v_{1}\bullet_{g}v_{2}-v_{2}\bullet_{g} v_{1}-2\pi \sqrt{-1}\res_{x_{0}}Y(e^{2\pi \sqrt{-1}x(a_{1}+\mathcal{N})}v_{1}, x_{0})v_{2}+\tilde{O}(V)=0.$$
Thus $\sum_{i=1}^{k}\rho(u^{i}_{1}+\tilde{O}(V))
\cdots \rho(u^{i}_{k}+\tilde{O}(V))w=0$. We have a ontracdiction. So $S_{g}^{1}(M)_{0}$ is linearly isomorphic to $M$. We shall 
identify $S_{g}^{1}(M)_{0}$ with $M$. Then $S_{g}^{1}(M)=S_{g}^{1}(M)_{+}\otimes M$.

Let $\mathcal{J}$ be the $\overline{\C}_{+}\times \C/\Z$-graded 
$T(L(V, g))$-submodule of $S_{g}^{1}(M)$
generated by the coefficients in $x$ of  
$$Y_{0}^{t}(L(-1)u, x)w-\frac{d}{dx}Y_{0}^{t}(u, x)w+x^{-1}Y_{0}^{t}(\mathcal{N}u, x)w$$
and the coefficients in $x_{0}$ and $x_{2}$ of 
$$(x_0+x_2)^{l} Y^{t}_{0}(u, x_0+x_2)Y^{t}_{0}(v, x_2)w
- (x_2+x_0)^{l} Y^{t}_{0}\left(Y\left(\left(1+\frac{x_0}{x_2}\right)^{\mathcal{N}}u,
x_0\right)v, x_2\right)w$$
(which are indeed in $S_{g}^{1}(M)$)
for $u\in V^{[a+\Z]}, v\in V$, $w\in S_{g}^{1}(M)$ and 
$l\in a+\Z$ satisfying $u(n)w=0$ for $n\in a+\Z$ such that $n-a\ge l-a$. 

Let $S_{g}(M)=S_{g}^{1}(M)/\mathcal{J}$. Then $S_{g}(M)$ is also 
a $\overline{\C}_{+}\times \C/\Z$-graded 
$T(L(V, g))$-module. We can
still use elements of $T(L(V, g))\otimes M$ 
to represent elements of $S_{g}(M)$. But note that these 
elements now satisfy relations. We equip $S_{g}(M)$ with 
the vertex operator map $Y_{0}^{g}: V\otimes S_{g}(M)\to S_{g}(M)[[x, x^{-1}]]$
given by $u\otimes w\mapsto Y_{0}^{g}(u, x)w=Y_{0}^{t}(u, x)w$. 
As in $S_{g}^{1}(M)$, for $u\in V^{[a+\Z]}$, $m\in a+\Z$ and $w\in S_{g}(M)$,
we also have $u(m)w=0$ when the real part of $\wt u-m-1$ is sufficiently negative. Clearly 
$Y_{0}^{g}(\mathbf{1}, x)=1_{S_{g}(M)}$ (where $1_{S_{g}(M)}$ is the identity 
operator on $S_{g}(M)$). 
By definition, we know that the commutator formula, the weak associativity 
and the $L(-1)$-derivative property for $Y_{0}^{g}$ all hold. The  $\overline{\C}_{+}$-grading condition 
is also clear.
Thus $S_{g}(M)$ is a  $\overline{\C}_{+}$-graded weak $g$-twisted $V$-module. 

We still need to show that $\Omega_{g}(S_{g}(M))=M$. To prove  this fact,
we need only prove that the relations given by $\mathcal{J}$ already hold in $M$.
Let $\pi_{M}$ be the projection from $S_{g}^{1}(M)$ to $M$. we need only prove that 
\begin{equation}\label{L(-1)-M}
\pi_{M}Y_{0}^{t}(L(-1)u, x)w-\pi_{M}\frac{d}{dx}Y_{0}^{t}(u, x)w+\pi_{M}x^{-1}Y_{0}^{t}(\mathcal{N}u, x)w=0
\end{equation}
and 
\begin{equation}\label{weak-assoc-M}
(x_0+x_2)^{l} \pi_{M}Y^{t}_{0}(u, x_0+x_2)Y^{t}_{0}(v, x_2)w
= (x_2+x_0)^{l} \pi_{M}Y^{t}_{0}\left(Y\left(\left(1+\frac{x_0}{x_2}\right)^{\mathcal{N}}u,
x_0\right)v, x_2\right)w
\end{equation}
hold in $S_{g}^{1}(M)$ for homogeneous  $u\in V^{[a+\Z]}, v\in V$, $w\in S_{g}^{1}(M)$.  Recall 
$l(u)\in a+\Z$ in Section 4 defined by 
$\wt u-n-1<0$ (thus $u(n)w=0$ for $w\in M$) for $n\in a+\Z$ satisfying $n-a\ge l(u)-a$. We first prove these formulas
for $w\in M\subset  S_{g}^{1}(M)$ and $l=l(u)$.

Since $u$ is homogeneous with respect to the $\C$-grading (weight grading),  the terms whose weights are in $\mathbb{I}$  in the series
$Y_{0}^{t}(L(-1)u, x)$ , $\frac{d}{dx}Y_{0}^{t}(u, x)$ and $x^{-1}Y_{0}^{t}(\mathcal{N}u, x)$
 are $0$ when $a\not \in P(V)\cap \mathbb{I}$ and are $(L(-1)u)(\wt u+a)x^{-\swt u-a-1}$, $(-\wt u-a) u(\wt u+a-1)x^{-\swt u-a-1}$ and 
$(\mathcal{N}u)(\wt u+a-1)x^{-\swt u-a-1}$, respectively, when $a\in P(V)\cap \mathbb{I}$. Then in $S_{g}^{1}(M)$,
for $a\in P(V)\cap \mathbb{I}$, 
\begin{eqnarray}\label{L(-1)-M-1}
\lefteqn{\pi_{M}Y_{0}^{t}(L(-1)u, x)w-\pi_{M}\frac{d}{dx}Y_{0}^{t}(u, x)w+\pi_{M}x^{-1}Y_{0}^{t}(\mathcal{N}u, x)w}\nn
&&=((L(-1)u)(\wt u+a)+(\wt u+a) u(\wt u+a-1)+(\mathcal{N}u)(\wt u+a-1))wx^{-\swt u-a-1}\nn
&&=o_{t}((L(-1)+L(0)+a+\mathcal{N})u)wx^{-\swt u-a-1}\nn
&&=o_{t}((L(-1)+L(0)+a+\mathcal{N})\mathcal{U}(1) (\mathcal{U}(1))^{-1}u)wx^{-\swt u-a-1}.
\end{eqnarray}
By (1.15) in \cite{H1}  ($(L(-1)+L(0))\mathcal{U}(1) =\mathcal{U}(1) \frac{L(-1)}{2\pi \sqrt{-1}}$)
and the fact that the Virasoro operators commutes with $\mathcal{N}$, 
we have 
$$(L(-1)+L(0)+a+\mathcal{N})\mathcal{U}(1) =\mathcal{U}(1) \left(\frac{L(-1)}{2\pi \sqrt{-1}}+a+\mathcal{N}\right).$$
Thus the right-hand side of (\ref{L(-1)-M-1}) is equal to 
\begin{eqnarray}\label{L(-1)-M-2}
\lefteqn{o_{t}\left(\mathcal{U}(1) \left(\frac{L(-1)}{2\pi \sqrt{-1}}+a+\mathcal{N}\right)(\mathcal{U}(1))^{-1}u\right)wx^{-\swt u-a-1}}\nn
&&=\rho\left(\left(\frac{L(-1)}{2\pi \sqrt{-1}}+a+\mathcal{N}\right)(\mathcal{U}(1))^{-1}u+\tilde{O}_{g}(V)\right)wx^{-\swt u-a-1}.
\end{eqnarray}
Since $\left(\frac{L(-1)}{2\pi \sqrt{-1}}+a+\mathcal{N}\right)v\in \tilde{O}_{g}(V)$ for $v\in V$ (by Theorem \ref{tilde-A-g}) and
$\rho(u+\tilde{O}_{g}(V))=0$ for $u\in \tilde{O}_{g}(V)$, the right-hand side of (\ref{L(-1)-M-2}) and thus the 
right-hand side of (\ref{L(-1)-M-1}) are $0$ in $S_{g}^{1}(M)$, proving (\ref{L(-1)-M}) in the case $w\in M$. 

For the same $u$ and $v$ but for $w\in S^{1}_{g}(V)_{+}$, 
by straightforward calculations using the commutator formula for $Y_{0}^{t}$, the properties of the formal $\delta$-function
and $\res$ and the $L(-1)$-derivative property for $Y$, it is easy to show that 
$$Y_{0}^{t}(L(-1)u, x)w-\frac{d}{dx}Y_{0}^{t}(u, x)w+x^{-1}Y_{0}^{t}(\mathcal{N}u, x)$$
commutes with $Y_{0}^{t}(v, y)$ for $v\in V$ and thus also commutes with $v(n)$ for $v\in V^{[b+\Z]}$, $b\in P(V)$ and 
$n\in b+\Z$. Here we omit the calculations. Now let $w=v_{1}(m_{1})\cdots v_{k}(m_{k})w_{0}\in S^{1}_{g}(M)_{+}$ for 
$v_{i}\in V^{[a_{i}+\Z]}$, $m_{i}\in a_{i}+\Z$, $w_{0}\in M$ such that $\wt v_{i}-m_{i}-1>0$.
Then 
\begin{eqnarray}\label{L(-1)-M-1}
\lefteqn{\pi_{M}\left(Y_{0}^{t}(L(-1)u, x)w-\frac{d}{dx}Y_{0}^{t}(u, x)w+x^{-1}Y_{0}^{t}(\mathcal{N}u, x)\right)w}\nn
&&=\pi_{M}\left(Y_{0}^{t}(L(-1)u, x)w-\frac{d}{dx}Y_{0}^{t}(u, x)w+x^{-1}Y_{0}^{t}(\mathcal{N}u, x)\right)
v_{1}(m_{1})\cdots v_{k}(m_{k})w_{0}\nn
&&=\pi_{M}v_{1}(m_{1})\cdots v_{k}(m_{k})\left(Y_{0}^{t}(L(-1)u, x)w-\frac{d}{dx}Y_{0}^{t}(u, x)w+x^{-1}Y_{0}^{t}(\mathcal{N}u, x)\right)
w_{0}.\nn
\end{eqnarray}
Since the weight of $v_{1}(m_{1})\cdots v_{k}(m_{k})$ is bigger than $0$ and $w_{0}\in M$, 
the right-hand side of (\ref{L(-1)-M-1}) must be $0$, proving (\ref{L(-1)-M}) in this case. 

We now prove (\ref{weak-assoc-M}) in the case $w\in M$. To prove (\ref{weak-assoc-M}) in this case, we need only prove that coefficients in 
$x_{0}$ of the two sides of (\ref{weak-assoc-M})
are equal, that is, for $n\in \Z$,
\begin{eqnarray}\label{weak-assoc-M-3}
\lefteqn{\res_{x_{0}}x_{0}^{n}(x_{0}+x_{2})^{l(u)}
\pi_{M}Y^{t}_{0}(u, x_0+x_2)Y^{t}_{0}(v, x_2)w}\nn
&&=\res_{x_{0}}x_{0}^{n}(x_{2}+x_{0})^{l(u)}
 \pi_{M}Y_{0}^{t}\left(Y\left(\left(1+\frac{x_{0}}{x_{2}}\right)^{\mathcal{N}}u, x_{0}\right)v, x_{2}\right)w.
\end{eqnarray}

For $n\in \Z$, by changing the variable from $x_{0}$ to $x_{1}$, we obtain
\begin{eqnarray}\label{weak-assoc-M-0}
\lefteqn{\res_{x_{0}}x_{0}^{n}(x_{0}+x_{2})^{l(u)}
\pi_{M}Y^{t}_{0}(u, x_0+x_2)Y^{t}_{0}(v, x_2)w}\nn
&&=\res_{x_{1}} (x_{1}-x_{2})^{n}x_{1}^{l(u)}\pi_{M}Y^{t}_{0}(u, x_1)Y^{t}_{0}(v, x_2)w.
\end{eqnarray}
 For $n\in \N$, using the commutator formula for $Y^{t}_{0}$ for
$S_{g}^{1}(M)$, we see that the right-hand side of (\ref{weak-assoc-M-0}) is equal to 
\begin{eqnarray}\label{weak-assoc-M-1}
\lefteqn{\res_{x_{1}} (x_{1}-x_{2})^{n}x_{1}^{l(u)}\pi_{M}Y^{t}_{0}(v, x_2)Y^{t}_{0}(u, x_1)w}\nn
&&\quad +\res_{x_{1}}
\res_{x_{0}}(x_{1}-x_{2})^{n}x_{1}^{l(u)}x_{1}^{-1}\delta\left(\frac{x_{2}+x_{0}}{x_{1}}\right)\left(\frac{x_{2}+x_{0}}{x_{1}}\right)^{a}\cdot\nn
&&\quad\quad\quad\quad\quad\quad\quad\cdot 
 \pi_{M}Y_{0}^{t}\left(Y\left(\left(1+\frac{x_{0}}{x_{2}}\right)^{\mathcal{N}}u, x_{0}\right)v, x_{2}\right)w.
\end{eqnarray}
Since $u(l(u))w=0$,
$$\res_{x_{1}}x_{1}^{l(u)+i}Y^{t}_{0}(u, x_1)w=0$$
for $0\le i\le n$.
Thus the first term in the right-hand side of (\ref{weak-assoc-M-1}) is $0$. 
The second term in the right-hand side of (\ref{weak-assoc-M-1}) is equal to 
\begin{equation}\label{weak-assoc-M-2}
\res_{x_{0}}x_{0}^{n}(x_{2}+x_{0})^{l(u)}
 \pi_{M}Y_{0}^{t}\left(Y\left(\left(1+\frac{x_{0}}{x_{2}}\right)^{\mathcal{N}}u, x_{0}\right)v, x_{2}\right)w.
\end{equation}
From (\ref{weak-assoc-M-1}) and (\ref{weak-assoc-M-2}), we obtain 
(\ref{weak-assoc-M-3}) in the case $n\in \N$. 

For $n\in -\Z_{+}-1$, 
$\res_{x_{1}}(x_{1}-x_{2})^{n}x_{1}^{l(u)}Y^{t}_{0}(u, x_1)$
contains terms proportional to $u(l(u)+n-i)$ for $i\in \N$. Recall from Section 4 that $l(u)=\wt u+a$ when $a\in \mathbb{I}$ and 
$l(u)=a+\wt u-1$ when $a\in P(V)\cap((0, 1)+\mathbb{I})$.  Together with $n\in -\Z_{+}-1$
and $i\in \N$, we have $\pi_{M}u(l(u)+n-i)v(k)w=0$ for $k\in \C$. Thus the right-hand side of 
(\ref{weak-assoc-M-0}) is $0$. We now show that the right-hand side of (\ref{weak-assoc-M-3}) 
is also $0$. First note that for homogeneous $u\in V^{[a+\Z]}$ and $w\in M$, 
$\pi_{M}Y_{0}^{g}(x^{L(0)}u, x)w=u(\wt u+a-1)wx^{-a}=o_{t}(u)wx^{-a}$ if $a\in P(V)\cap \mathbb{I}$ and is $0$ if 
otherwise. Now come back to $u, v, w$ as above.
Recall the operator $\L$ introduced in Section 3. 
By changing the variable $x_{0}$ to $y=\frac{x_{0}}{x_{2}}$, we obtain
\begin{eqnarray}\label{weak-assoc-M-4}
\lefteqn{\res_{x_{0}}x_{0}^{n}(x_{2}+x_{0})^{l(u)}
 \pi_{M}Y_{0}^{t}\left(Y\left(\left(1+\frac{x_{0}}{x_{2}}\right)^{\mathcal{N}}u, x_{0}\right)v, x_{2}\right)w}\nn
&&=\res_{y}x_{2}(x_{2}y)^{n}(x_{2}+x_{2}y)^{l(u)}
 \pi_{M}Y_{0}^{t}\left(Y\left((1+y)^{\mathcal{N}}u, x_{2}y\right)v, x_{2}\right)w\nn
&&=x_{2}^{l(u)+n+1+\swt u+\swt v}
 \pi_{M}Y_{0}^{t}\left(x_{2}^{L(0)}\res_{y}y^{n}Y\left((1+y)^{\L+\mathcal{N}}u, y\right)v, x_{2}\right)w\nn
&&=x_{2}^{l(u)+n+1+\swt u+\swt v}
 o_{t}(\res_{y}y^{n}Y\left((1+y)^{\L+\mathcal{N}}u, y\right)v)w\nn
&&=x_{2}^{l(u)+n+1+\swt u+\swt v}
\rho(\mathcal{U}(1)^{-1}\res_{y}y^{n}Y\left((1+y)^{\L+\mathcal{N}}u, y\right)v+\tilde{O}_{g}(V))w.
\end{eqnarray}
By (\ref{bullet=*}) and (\ref{tilde-O-g}), 
\begin{eqnarray}\label{weak-assoc-M-4.1}
\lefteqn{\mathcal{U}(1)^{-1}\res_{y}y^{n}Y\left((1+y)^{\L+\mathcal{N}}u, y\right)v}\nn
&&=\res_{y}y^{n}Y\left((1+y)^{\L-L(0)+\mathcal{N}}
\mathcal{U}(1)^{-1}u, \frac{1}{2\pi \sqrt{-1}}\log (1+y)\right)\mathcal{U}(1)^{-1}v\nn
&&\in \tilde{O}_{g}(V).
\end{eqnarray}
Thus the right-hand side of (\ref{weak-assoc-M-4}) is $0$. 

We still need to prove (\ref{weak-assoc-M-3}) in the case $n=-1$. We first prove this formula when $a\in P(V)\cap ((0, 1)+\mathbb{I})$.
Since $l(u)=\wt u+a-1$ and $\pi_{M}u(l(u)-1-i)v(k)w=\pi_{M}u(\wt u+a-2-i)v(k)w=0$ for $i\in \N$ and $k\in \C$,
\begin{eqnarray}\label{weak-assoc-M-5}
\lefteqn{\res_{x_{0}}x_{0}^{-1}(x_{0}+x_{2})^{l(u)}
\pi_{M}Y^{t}_{0}(u, x_0+x_2)Y^{t}_{0}(v, x_2)w}\nn
&&=\res_{x_{1}} (x_{1}-x_{2})^{-1}x_{1}^{l(u)}\pi_{M}Y^{t}_{0}(u, x_1)Y^{t}_{0}(v, x_2)w\nn
&&=0.
\end{eqnarray}
On the other hand, the same calculation as in (\ref{weak-assoc-M-4}) and (\ref{weak-assoc-M-4.1}) gives
\begin{eqnarray*}
\lefteqn{\res_{x_{0}}x_{0}^{-1}(x_{2}+x_{0})^{l(u)}
 \pi_{M}Y_{0}^{t}\left(Y\left(\left(1+\frac{x_{0}}{x_{2}}\right)^{\mathcal{N}}u, x_{0}\right)v, x_{2}\right)w}\nn
&&=x_{2}^{l(u)+\swt u+\swt v}
\rho(\mathcal{U}(1)^{-1}\res_{y}y^{-1}Y\left((1+y)^{\L+\mathcal{N}}u, y\right)v+\tilde{O}_{g}(V))w\nn
\end{eqnarray*}
\begin{eqnarray}\label{weak-assoc-M-6}
&&=x_{2}^{l(u)+\swt u+\swt v}
\rho\Biggl(\res_{y}y^{-2}Y\Biggl(y^{L(0)-\L+\A}(1+y)^{\L-L(0)+\mathcal{N}}\cdot\nn
&&\quad\quad\quad\quad\quad\quad\quad\quad\quad\quad\quad\quad\cdot
\mathcal{U}(1)^{-1}u, \frac{1}{2\pi \sqrt{-1}}\log (1+y)\Biggr)\mathcal{U}(1)^{-1}v+\tilde{O}_{g}(V)\Biggr)w\nn
&&=0
\end{eqnarray}
because $a\in P(V)\cap ((0, 1)+\mathbb{I})$, proving (\ref{weak-assoc-M-3}) in this case ($n=-1$ and $a\in P(V)\cap ((0, 1)+\mathbb{I})$).

In the case $a\in P(V)\cap \mathbb{I}$, since $l(u)=\wt u+a$, we have 
$\pi_{M}u(l(u)-1-i)v(k)w=\pi_{M}u(\wt u+a-1-i)v(k)w=0$ for $i\in \Z_{+}$ and  $k\in \C$,
$\pi_{M}u(l(u)-1)v(k)w=0$ for either $v\in V^{[b]}$, $b\in P(V)\cap ((0, 1)+\mathbb{I})$, $k\in \C$ or 
$v\in V^{[b]}$, $b\in P(V)\cap \mathbb{I}$, $k\ne \wt v+b-1$
and 
$\pi_{M}u(l(u)-1)v(\wt v+b-1)w=u(\wt u+a-1)v(\wt v+b-1)w$ for $v\in V^{[b]}$, $b\in P(V)\cap \mathbb{I}$.
Using these formulas, we see that  in the case $v\in V^{[b]}$ and 
$b\in P(V)\cap ((0, 1)+\mathbb{I})$, the right-hand side of  (\ref{weak-assoc-M-0}) is equal to 
$0$. But in this case, 
$$Y^{t}_{0}\left(Y\left(\left(1+\frac{x_{0}}{x_{2}}\right)^{\mathcal{N}}u, x_{0}\right)v, x_{2}\right)w=0.$$
Thus (\ref{weak-assoc-M-3}) holds in this case ($n=-1$, $a\in P(V)\cap \mathbb{I}$ and $b\in P(V)\cap ((0, 1)+\mathbb{I})$). 
In the case $v\in V^{[b]}$ and $b\in P(V)\cap \mathbb{I}$, the right-hand side of  (\ref{weak-assoc-M-0}) is equal to 
\begin{eqnarray*}
\lefteqn{x_{2}^{-\swt v} u(\wt u+a-1)v(\wt v+b-1)w}\nn
&&= x_{2}^{-\swt v}o_{t}(u)o_{t}(v)w\nn
&&= x_{2}^{-\swt v}\rho(\mathcal{U}(1)^{-1}u+\tilde{O}_{g}(V))\rho(\rho(\mathcal{U}(1)^{-1}v+\tilde{O}_{g}(V))w \nn
&&= x_{2}^{-\swt v}\rho((\mathcal{U}(1)^{-1}u\bullet_{g} \mathcal{U}(1)^{-1}v)+\tilde{O}_{g}(V))w \nn
&&=x_{2}^{-\swt v}\rho(\mathcal{U}(1)^{-1}(u*_{g} v)+\tilde{O}_{g}(V))w\nn
&&= x_{2}^{-\swt v} o_{t}(u*_{g} v)w\nn
&&=x_{2}^{-\swt v}\pi_{M}Y^{t}_{0}(x_{2}^{L(0)}(u*_{g} v), x_{2})w
\end{eqnarray*}
\begin{eqnarray*}
&&= x_{2}^{-\swt v}\res_{y}y^{-1}\pi_{M}Y^{t}_{0}(x_{2}^{L(0)}Y((1+y)^{L(0)+\mathcal{N}}u, y)v, x_{2})w\nn
&&=\res_{x_{0}}x_{0}^{-1} x_{2}^{-\swt v}\pi_{M}Y^{t}_{0}\left(x_{2}^{L(0)}
Y\left(\left(1+\frac{x_{0}}{x_{2}}\right)^{L(0)+\mathcal{N}}u, \frac{x_{0}}{x_{2}}\right)v, x_{2}\right)w\nn
&&=\res_{x_{0}}x_{0}^{-1} (x_{2}+x_{0})^{l}\pi_{M}Y^{t}_{0}\left(
Y\left(\left(1+\frac{x_{0}}{x_{2}}\right)^{\mathcal{N}}u, x_{0}\right)v, x_{2}\right)w,
\end{eqnarray*}
proving  (\ref{weak-assoc-M-3}) in the case $n=-1$, $a, b\in P(V)\cap \mathbb{I}$.

We now prove (\ref{weak-assoc-M}) 
for the same $u$ and $v$ but for $w\in S^{1}_{g}(V)_{+}$. By straightforward calculations generalizing
the last part of the calculations in the proof of Proposition 6.1 in \cite{DLM}, for $u\in V^{[a+\Z]}$,
$a\in P(V)$, $v\in V$, $v_{1}\in V^{[a_{1}+\Z]}$, $a_{1}\in P(V)$, $m_{1}\in a_{1}+\Z$ such that 
$\Re(\wt v_{1}-m_{1}-1)>0$, 
$w\in S^{1}_{g}(V)$, $l\in a+\Z$ such that 
$u(n)v_{1}(m_{1})w=0$ for $n\in a+\Z$ and $n-a\ge l-a$, we have
\begin{eqnarray}\label{weak-assoc-diff}
\lefteqn{\Biggl((x_0+x_2)^{l} Y^{t}_{0}(u, x_0+x_2)Y^{t}_{0}(v, x_2)}\nn
&&\quad \quad - (x_2+x_0)^{l} Y^{t}_{0}\left(Y\left(\left(1+\frac{x_0}{x_2}\right)^{\mathcal{N}}u,
x_0\right)v, x_2\right)\Biggr)v_{1}(m_{1})w\nn
&&=v_{1}(m_{1})\Biggl((x_0+x_2)^{l} Y^{t}_{0}(u, x_0+x_2)Y^{t}_{0}(v, x_2)\nn
&&\quad\quad\quad\quad\quad\quad\quad\quad - (x_2+x_0)^{l} Y^{t}_{0}\left(Y\left(\left(1+\frac{x_0}{x_2}\right)^{\mathcal{N}}u,
x_0\right)v, x_2\right)\Biggr)w
\end{eqnarray}
Here we omit the detailed calculations. Now let $w=v_{1}(m_{1})\cdots v_{k}(m_{k})w_{0}$
where $v_{i}\in V^{[a_{i}+\Z]}$, $m_{i}\in a_{i}+\Z$ satisfying $\Re(\wt v_{i}-m_{i}-1)>0$
and $w_{0}\in M$.  Then 
by (\ref{weak-assoc-diff}), for $l\in a+\Z$ such that $u(n)v_{1}(m_{1})\cdots v_{k}(m_{k})w_{0}=0$ for $n\in a+\Z$ and $n-a\ge l-a$,
we have
\begin{eqnarray*}
\lefteqn{\pi_{M}\Biggl((x_0+x_2)^{l} Y^{t}_{0}(u, x_0+x_2)Y^{t}_{0}(v, x_2)}\nn
&&\quad \quad - (x_2+x_0)^{l} Y^{t}_{0}\left(Y\left(\left(1+\frac{x_0}{x_2}\right)^{\mathcal{N}}u,
x_0\right)v, x_2\right)\Biggr)v_{1}(m_{1})\cdots v_{k}(m_{k})w_{0}\nn
&&=\pi_{M}v_{1}(m_{1})\cdots v_{k}(m_{k})\Biggl((x_0+x_2)^{l} Y^{t}_{0}(u, x_0+x_2)Y^{t}_{0}(v, x_2)\nn
&&\quad\quad\quad\quad\quad\quad\quad\quad - (x_2+x_0)^{l} Y^{t}_{0}\left(Y\left(\left(1+\frac{x_0}{x_2}\right)^{\mathcal{N}}u,
x_0\right)v, x_2\right)\Biggr)w_{0}\nn
&&=0,
\end{eqnarray*}
where the last step uses the fact that $\Re(\wt v_{i}-m_{i}-1)>0$ for $i=1, \dots, k$ 
together with the projection $\pi_{M}$ force
us to take only components of 
$$(x_0+x_2)^{l} Y^{t}_{0}(u, x_0+x_2)Y^{t}_{0}(v, x_2)- (x_2+x_0)^{l} Y^{t}_{0}\left(Y\left(\left(1+\frac{x_0}{x_2}\right)^{\mathcal{N}}u,
x_0\right)v, x_2\right)$$
of negative weights. This proves (\ref{weak-assoc-M})  in the case that $w\in S^{1}_{g}(V)_{+}$.

Let $W$ be a $\overline{\C}_{+}$-graded weak $g$-twisted $V$-module. 
We define a linear map from $S_{g}(\Omega_{g}(W))$ to $W$ by 
mapping $u_{1}(m_{1})\cdots u_{k}(m_{k})w$ of $S_{g}(M)$
to $u_{1}(m_{1})\cdots u_{k}(m_{k})w$  of $W$ for 
$u_{i}\in V^{[a_{i}+\Z]}$, $m_{i}\in a_{i}+\mathbb{Z}$ ($i=1, \dots, k$) and $w\in \Omega_{g}(W)$.
Note that the only relations among $u_{1}(m_{1})\cdots u_{k}(m_{k})w$
for $u_{i}\in V$, $m_{i}\in \mathbb{Z}$ ($i=1, \dots, k$) and $w\in \Omega_{g}(W)$
are given by the action of $V$ on $\Omega_{g}(W)$, the commutator formula,
the weak associativity and the $L(-1)$-derivative property 
for $Y_{0}^{g}$. These relations also hold in $W$. Thus 
this map is well-defined. Clearly, this is a surjective homomorphism of 
$\overline{\C}_{+}$-graded weak $g$-twisted $V$-modules from $S_{g}(\Omega_{g}(W))$ to 
$\overline{\C}_{+}$-graded weak $g$-twisted $V$-submodule of $W$
generated by $\Omega_{g}(W)$.

The other statements follow immediately.
\epfv

\section{Isomorphisms between $Z_{g}(V)$, $\tilde{A}_{g}(V)$ and $A_{g}(V)$}

In this section, using the results obtained in the preceding sections, 
we prove that the algebras $Z_{g}(V)$, $\tilde{A}_{g}(V)$ and $A_{g}(V)$ are isomorphic to each other. 

\begin{thm}\label{isomorphism}
The associative algebras $Z_{g}(V)$ and $A_{g}(V)$ are isomorphic.
\end{thm}
\pf 
We view $A_{g}(V)$ as an $A_{g}(V)$-module. By Theorem 
\ref{functor-f}, $S_{g}(A_{g}(V))$ is a $\overline{\C}_{+}$-graded $g$-twisted $V$-module such that 
$\Omega_{g}(S_{g}(A_{g}(V)))=A_{g}(V)$. Then by Proposition \ref{z-module},
$A_{g}(V)$ is a $Z_{g}(V)$-module and, by definition, the action of $Z_{g}(V)$ on 
$A_{g}(V)$ is given by the homomorphism from $Z_{g}(V)$ to the endomorphism 
algebra of $A_{g}(V)$ determined by $[o_{g}(v)]\mapsto v+O_{g}(V)$ for $v\in V$, where 
$v+O_{g}(V)$ is viewed as an element of the endomorphism algebra of $A_{g}(V)$.
Here we have also used  Proposition \ref{elements-z-g}. 
Since the image of this homomorphsim is in fact in $A_{g}(V)$, it is a homomorphism
$f: Z_{g}(V)\to A_{g}(V)$ of associative algebras. 

Conversely, we view $Z_{g}(V)$ as a $Z_{g}(V)$-module. By Theorem \ref{functor-H},
$H_{g}(Z_{g}(V))$ is a $\overline{\C}_{+}$-graded $g$-twisted $V$-module such that 
$\Omega_{g}(H_{g}(Z_{g}(V)))=Z_{g}(V)$. Then by Theorem \ref{module},
$Z_{g}(V)$ is an $A_{g}(V)$-module and, by definition, the action of $A_{g}(V)$ on 
$Z_{g}(V)$ is given by the homomorphism from $A_{g}(V)$ to the endomorphism 
algebra of $Z_{g}(V)$ defined by $v+O_{g}(V)\mapsto [o_{g}(v)]$ for $v\in V$, where 
$[o_{g}(v)]$ is viewed as an element of the endomorphism algebra of $Z_{g}(V)$.
Since the image of this homomorphsim is in fact in $Z_{g}(V)$, it is a homomorphism
$g: A_{g}(V)\to Z_{g}(V)$ of associative algebras. 

Clearly $f$ and $g$ are inverse to each other. Thus they are isomorphisms of associative algebras.
\epfv

\begin{cor}
The associative algebras $Z_{g}(V)$ and $\tilde{A}_{g}(V)$ are isomorphic.\epf
\end{cor}

From Theorem \ref{isomorphism}, we also obtain an improvement of Proposition \ref{elements-z-g}:

\begin{cor}\label{elements-z-g-1}
The elements of $Z_{g}(V)$ are of the form 
$[o_{g}(v)]$ for $v \in V$. 
\end{cor}
\pf
Since an element of $A_{g}(V)$ is of the form $v+O_{g}(V)$ for $v \in V$ and under the isomorphism from 
$A_{g}(V)$ to $Z_{g}(V)$, such an element is mapped to $[o_{g}(v)]$ for $v \in V$, 
the conclusion follows.
\epfv

\begin{rema}
{\rm Because of Theorem \ref{isomorphism}, if we have results for the algebra $Z_{g}(V)$
(or $A_{g}(V)$), we can 
immediately conclude that the same results hold for $A_{g}(V)$ (or $Z_{g}(V)$). 
In particular, we can derive the results in Section 5 on $A_{g}(V)$ (or $Z_{g}(V)$) 
from the corresponding results on $Z_{g}(V)$ (or $A_{g}(V)$). But since our proof of 
Theorem \ref{isomorphism} has used the results in Section 5 for both $Z_{g}(V)$ and $A_{g}(V)$,
the proofs of the results in Section 5 are all needed. }
\end{rema}

\begin{rema}
{\rm In many cases, using Corollary \ref{elements-z-g-1},
 it is easier to calculate $Z_{g}(V)$ than to calculate $A_{g}(V)$ or $\tilde{A}_{g}(V)$. The 
algebra $\tilde{A}_{g}(V)$ is more natural for the study of modular invariance. The algebra $A_{g}(V)$ 
is in some sense a bridge between $Z_{g}(V)$ and $\tilde{A}_{g}(V)$.
We expect that $Z_{g}(V)$ will play an important role in the future construction 
and study of twisted modules and twisted intertwining operators. }
\end{rema}

\noindent {\small \sc Department of Mathematics, Rutgers University,
110 Frelinghuysen Rd., Piscataway, NJ 08854-8019}

\noindent {\em E-mail address}: {\tt yzhuang@math.rutgers.edu}

\vspace{1em}

\noindent {\small \sc Department of Mathematics, University of Notre Dame,
278 Hurley Building, Notre Dame, IN 46556}

\noindent {\em E-mail address}:  {\tt  jyang7@nd.edu}

\end{document}